\documentclass[reqno,11pt]{amsart}
\usepackage[utf8]{inputenc}

\usepackage{amsmath,amsxtra,amsbsy,enumerate, latexsym, amsfonts, amssymb, amsthm, amscd, stmaryrd}
\usepackage{graphics,epsf,psfrag}
\usepackage{bbm, dsfont}
\usepackage{mathtools}
\usepackage{xcolor}
\usepackage{float}
\usepackage{tikz}
\usepackage[headings]{fullpage}
\usepackage{colonequals}

\numberwithin{equation}{section}
\allowdisplaybreaks

\usepackage{xcolor}
\usepackage[normalem]{ulem}
\usepackage{hyperref}

\hypersetup{
    colorlinks=true,
    linkcolor=blue,
    filecolor=magenta,      
    urlcolor=blue,
    citecolor=blue,
     pdfstartview=FitV,
    pdfpagemode=FullScreen,
    hyperindex,breaklinks
    }

\makeatletter
\DeclareUrlCommand\ULurl@@{%
  \def\UrlLeft{\uline\bgroup}%
  \def\UrlRight{\egroup}}
\def\ULurl@#1{\hyper@linkurl{\ULurl@@{#1}}{#1}}
\DeclareRobustCommand*\ULurl{\hyper@normalise\ULurl@}
\makeatother

\definecolor{darkred}{rgb}{0.7,0.1,0.1}

\definecolor{darkgreen}{rgb}{0.1,0.7,0.1}

\let\oldtocsection=\tocsection
\let\oldtocsubsection=\tocsubsection
\let\oldtocsubsubsection=\tocsubsubsection
\renewcommand{\tocsection}[2]{\hspace{0em}\oldtocsection{#1}{#2}}
\renewcommand{\tocsubsection}[2]{\hspace{1em}\oldtocsubsection{#1}{#2}}
\renewcommand{\tocsubsubsection}[2]{\hspace{2em}\oldtocsubsubsection{#1}{#2}}
\DeclareRobustCommand{\SkipTocEntry}[5]{}

\newcommand{\bbE}{{\ensuremath{\mathbb E}} }

\newcommand{\bbN}{{\ensuremath{\mathbb N}} }

\newcommand{\bbP}{{\ensuremath{\mathbb P}} }

\newcommand{\bbR}{{\ensuremath{\mathbb R}} }

\newcommand{\bbZ}{{\ensuremath{\mathbb Z}} }

\renewcommand{\epsilon}{\varepsilon}

\newcommand{\ga}{\alpha}

\newcommand{\gga}{\gamma}            

\newcommand{\gep}{\varepsilon}       

\newcommand{\gO}{\Omega}
\newcommand{\gl}{\lambda}

\newcommand{\ind}{\mathbf{1}}

\newcommand{\lint}{\llbracket}
\newcommand{\rint}{\rrbracket}

\newtheorem{theorem}{Theorem}[section]
\newtheorem{lemma}[theorem]{Lemma}
\newtheorem{proposition}[theorem]{Proposition}

\newtheorem{rem}[theorem]{Remark}

\newcommand{\RN}[1]{%
  \textup{\uppercase\expandafter{\romannumeral#1}}%
}

\newcommand{\Var}{\mathrm{Var}}

\newcommand{\Gap}{\mathrm{gap}}

\newcommand{\dd}{\mathrm{d}}

\renewcommand{\tilde}{\widetilde}
\renewcommand{\hat}{\widehat}

\makeatletter
\def\captionfont@{\footnotesize}
\def\captionheadfont@{\scshape}

\long\def\@makecaption#1#2{%
  \vspace{2mm}
  \setbox\@tempboxa\vbox{\color@setgroup
    \advance\hsize-6pc\noindent
    \captionfont@\captionheadfont@#1\@xp\@ifnotempty\@xp
        {\@cdr#2\@nil}{.\captionfont@\upshape\enspace#2}%
    \unskip\kern-6pc\par
    \global\setbox\@ne\lastbox\color@endgroup}%
  \ifhbox\@ne 
    \setbox\@ne\hbox{\unhbox\@ne\unskip\unskip\unpenalty\unkern}%
  \fi
  \ifdim\wd\@tempboxa=\z@ 
    \setbox\@ne\hbox to\columnwidth{\hss\kern-6pc\box\@ne\hss}%
  \else 
    \setbox\@ne\vbox{\unvbox\@tempboxa\parskip\z@skip
        \noindent\unhbox\@ne\advance\hsize-6pc\par}%
\fi
  \ifnum\@tempcnta<64 
    \addvspace\abovecaptionskip
    \moveright 3pc\box\@ne
  \else 
    \moveright 3pc\box\@ne
    \nobreak
    \vskip\belowcaptionskip
  \fi
\relax
}
\makeatother

\newcommand{\gap}{\mathrm{gap}}

\def\writefig#1 #2 #3 {\rlap{\kern #1 truecm
\raise #2 truecm \hbox{#3}}}

\title[Spectral gap \& principal eigenfunction  of RCM in a line segment]{The spectral gap and principal eigenfunction  of the random conductance model in a line segment}

\author[Shangjie Yang]{Shangjie Yang} \address{Shangjie Yang \hfill\break
{Instituto de Matemática e Estatística,   Universidade Federal Fluminese}
\hfill\break
Rua Prof. Marcos Waldemar de Freitas Reis, s/n,  Niterói, RJ,  CEP 24210-201,
Brasil.
}\email{syang@id.uff.br}

\begin{document}

\maketitle

\begin{abstract} In this paper, we study the spectral gap and principal eigenfunction of the random walk
 in the line segment $\lint 1, N \rint$ with conductances $c^{(N)}(x, x+1)_{1\le x<N}$ where $c^{(N)}(x, x+1)>0$ is the rate of the random walk jumping from  site $x$ to site $x+1$ and 
vice versa. Let $r^{(N)}(x,  x+1) \colonequals 1/c^{(N)}(x,  x+1)$ be the resistances, and under the assumption
\begin{equation*}
 \limsup_{N\to \infty}\, \frac{1}{N}\sup_{1< m \le N}\, \left|  \sum_{x=2}^m r^{(N)}(x-1, x)- (m-1) \right|\;=\;0\,,
\end{equation*} 
  we prove that the spectral gap, denoted by $\mathrm{gap}_{N}$, of the process satisfies $\mathrm{gap}_{N}=(1+o(1))\pi^2/N^2$ and the principal eigenfunction $g_N$ with $g_N(1)=1$ corresponding to the spectral gap  is well approximated by $h_N(x) := \cos\left( (x-1/2)\pi/N \right)$.
\end{abstract}


\section{Introduction}

\subsection{Model}
\begin{figure}[h]
 \centering
   \begin{tikzpicture}[scale=.4,font=\tiny]
     \draw[thick] (25,-1) -- (53,-1);

     \foreach \x in { 25, 29,...,53} {\draw (\x,-1.3) -- (\x,-1);}

     \node[below] at (25,-1.3) {$1$};
     \node[below] at (53,-1.3) {$N$};

       \draw  (40.8, -0.7)  edge[bend right,<->](37.2, -0.7);
       \node[below] at (38.7, 0.9) {$c(x-1,x)$};

       \draw[fill] (41,-1) circle [radius=0.2];
       \node[below] at (41, -1.2) {$x$};

    \draw  (44.8, -0.7)  edge[bend right, <->](41.2,-0.7);

    \node[below] at (43, 0.9) {$c(x,  x+1)$};

        \draw  (40.8, -0.7)  edge[bend right,->](37.2, -0.7);
       \node[below] at (38.7, 0.9) {$c(x-1, x)$};
 \end{tikzpicture}
 \caption{A graphical explanation for the random walk in conductances $(c(x, x+1))_{1 \le x<N}$:  at edge $\{x, x+1 \}$ there is a Poisson clock with rate $c(x, x+1)>0$ for all $1\le x<N$, and we swap the contents of the two sites $x, \, x+1$ when the  clock on the edge $\{x, x+1\}$ rings.\label{fig:SEPrandcond}  
 }
\end{figure}
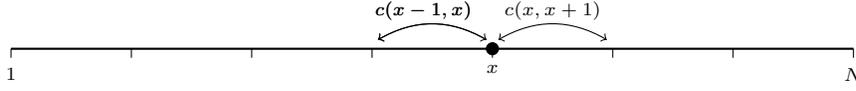

\noindent Let $\left(c^{(N)}(i, i+1)\right)_{1 \le i<N}$ be a sequence of strictly positive numbers. For $x<y$ real numbers, we define  
$\lint x,y\rint  \colonequals [x, y]\cap \bbZ$.
 For $N\in \bbN$, we consider 
the continuous-time random walk restricted in the segment  $\gO_N \colonequals \lint 1, N\rint$ 
with  its generator defined by ($f: \gO_{N} \mapsto \bbR$)
\begin{equation}\label{randwalk:generator}
(\Delta^{(c)} f)(x)\;:=\;  (c \nabla f)(x+1)- (c \nabla f)(x)\,, \quad \forall\, x \in \lint 1,\, N\rint
\end{equation}
where  we set
\begin{equation*}
\begin{gathered}
(c \nabla f)(x) \;\colonequals\; c^{(N)}(x-1,x) \left[f(x)-f(x-1) \right]\,, \quad x \in \lint 1, N+1 \rint\,,\\
c^{(N)}(N, N+1)[f(N+1)-f(N)]\;=\; c^{(N)}(0, 1)[f(1)-f(0)]\; \colonequals \;0\,.
\end{gathered}
\end{equation*}
This corresponds to Neumann boundary conditions at both endpoints. Throughout the paper, we often drop the superscript
 ``$(N)$" when it is clear in the context. 
We refer to Figure \ref{fig:SEPrandcond} for a graphical explanation. 
The model is a one-dimensional instance of the Random Conductance Model (RCM), first introduced in
 \cite{chernov1967replication}, with $(c^{(N)}(x, x+1))_{1\le x<N}$ being the conductances. 
The uniform probability measure 
$\mu_{N}$
on $\gO_{N}$,  defined by $\mu_{N}(x)=1/ N$ for all $x \in \gO_N$,  satisfies the detailed balance condition w.r.t.\  $\Delta^{(c)}$, and thus it is the unique invariant probability measure. We refer to \cite{Biskup2011survey, Kumagai2010SaintFlour}  and references therein for comprehensive introduction.

The spectral gap, denoted by $ \gap_{N}$,  of the generator $ \Delta^{(c)}$  which is the minimal strictly positive eigenvalue of $ -\Delta^{(c)}$.
To characterize it, we define the associated Dirichlet form  by ($f,g: \gO_{N} \mapsto \bbR$) 
\begin{equation*}
\mathcal{E}_{N}(f)\; \colonequals\; -\langle f, \Delta^{(c)} f\rangle_{\mu_{N}}\;=\;\sum_{x=1}^{N-1} \mu_{N}(x)c^{(N)}(x, x+1) \left[f(x+1)-f(x) \right]^2\,,
\end{equation*}
where $\langle f, g \rangle_{{\mu_{N}}} \colonequals \sum_{x \in \Omega_N} \mu_{N}(x)f(x)g(x)$ is the usual inner product in $L^2(\gO_{N}, \mu_{N})$. Moreover, the spectral gap $\gap_{N}$ is given by 
\begin{equation}\label{def:gap}
\gap_{N}\; \colonequals\;\,\, \inf_{f \ : \ \Var_{\mu_{N}}(f)>0}\,\,\frac{\mathcal{E}_{N}(f)}{\Var_{\mu_{N}}(f) }
\end{equation}  
where $\Var_{\mu_{N}}(f)\colonequals \langle f, f\rangle_{\mu_{N}}-\langle f, \ind \rangle^2_{\mu_{N}}$. The spectral gap is the rate of convergence to equilibrium (cf.\ \cite[Theorem 3.4]{Yang2021thesis}) and is also the best constant in the Poincar\'{e} inequality as shown in \eqref{def:gap}. The principal eigenfunction, corresponding to the spectral gap, can be exploited to obtain sharp bounds on the mixing time for the simple exclusion process and the  interchange process, see \cite{lacoin2016mixing, wilson2004mixing, yang2024cutoffSEP}.

 In the sequel, we always set
\begin{equation*}
 r^{(N)}(x, x+1)\;\colonequals\; 1/c^{(N)}(x, x+1),\quad  \forall\, x \in \lint 1, N-1\rint\,,
\end{equation*}
 which is referred to as resistance. For simplicity of notations, we drop the superscript ``$(N)$" when there is no confusion in the context.
 In the homogenous conductance, that is,   when $c(x, x+1) \equiv 1$ for all $x$, 
we know that 
\begin{equation*}
h_{i}^{(N)}(x)\;\colonequals\; \cos\left( \frac{ i \pi (x-1/2)}{N}\right)\,, \quad \forall\, x\, \in \lint 1, \,N \rint\,, \quad  \forall\,  i\in  \lint 0, N-1\rint
\end{equation*}
 are all the  eigenfunctions with corresponding eigenvalues $-2 (1-\cos(i \pi/N))$ for the generator $\Delta^{(c)}$  and thus the spectral is $\Gap_{N}=2 (1-\cos(\pi/N))=(1+o(1))\pi^2/N^2$.  A natural question is that  allowing arbitrary conductances with $c(x, x+1)>0$ for all $x$ rather than homogeneous conductances, how are the spectral gap and the principal eigenfunction affected by the \textit{disorder}?


\subsection{Monotonicity of the eigenfunctions} 
Given a function $f: \lint 1,N\rint\to \bbR$ and $2\le b\le c \le N-1$,
we say that $f$ admits a local maximum (respectively minimum) at $\lint b,c\rint$ if $f$ is constant on the interval $\lint b,c\rint$ and
$f(b-1)<f(b)$ and $f(c)> f(c+1)$ (resp.\ $f(b-1)>f(b)$ and $f(c)< f(c+1)$).
For  $j\ge 2$, we say that $f$ is 
$j$-monotone
 if  it displays exactly $(j-1)$ distinct local extrema in $\lint 2,N-1\rint$.

\begin{proposition}\label{prop:eigfunc}

For all conductance sequence $(c(x,x+1))_{1\le x<N}$ with $c(x,x+1)>0$, the following results hold:
\begin{enumerate}[(1)]

		\item\label{prop:egnfun1}
		 The operator  $-\Delta^{(c)}$ admits $N$ distinct eigenvalues: 
		\begin{equation}\label{eigv:distinct}
		0\;=\;\gl_0\;<\;\gl_1\;<\;\dots\;<\;\gl_{N-1}\,.
		\end{equation}
	\item 	If $(g_j)^{N-1}_{j=0}$ is a base of eigenfunctions of $-\Delta^{(c)}$ with respective eigenvalues $\gl_j$,
		then $g_1$ is  strictly monotone, and for $j\ge 2$, $g_j$ is $j$-monotone.
	\end{enumerate}
\end{proposition}

Note that $\Gap_N=\gl_1$.  Except that $g_1$ is strictly monotone,  Proposition \ref{prop:eigfunc} is a corollary of  \cite{Miclo2008eigenfunctions} which is concerned with the birth-and-death chain, and also provides nodal domain properties of eigenfunctions.
 We provide a self-contained proof for Proposition \ref{prop:eigfunc}, which sheds light on a sharp estimate for the spectral gap, used in the proof of Theorem \ref{th:gapshapeder}.
We comment that the strict monotonicity of $g_1$ can be exploited to prove the spectral gaps of the simple exclusion process and the interchange process are also $\gl_1$, see \cite{yang2024cutoffSEP}.

\subsection{Our result} 
Setting $r^{(N)}(n,m) \colonequals \sum^{m-1}_{x=n}r^{(N)}(x, x+1)$, 
 let us assume that
the sequence $(r^{(N)}(x, x+1))_{1\le x < N}$  satisfies 
\begin{equation}\label{LLN}
 \limsup_{N\to \infty} \delta_N^{(0)} \;\colonequals\;  \limsup_{N\to \infty}\, \frac{1}{N}\sup_{1\le n < m \le N} \left| (r^{(N)}(n,m)- (m-n) \right|\;=\;0\,.
\end{equation}

\begin{rem}
Note that the assumption \eqref{LLN} is equivalent to
\begin{equation}\label{eq32}
 \limsup_{N\to \infty}\, \frac{1}{N}\sup_{2\le m  \le N}\, \left| (r^{(N)}(1,m)- (m-1) \right|\;=\;0\,.
\end{equation}
When the sequence $(r^{(N)}(x, x+1))_{ 1 \le x < N}$ is  pairwise independent identically distributed random variables
 with common law denoted by $\bbP$ and its expectation $\bbE[r(x, x+1)]=1$, by \cite[Theorem 2.4.1]{Durrett} we have  
\begin{equation*}
\bbP \left( \lim_{N \to \infty}\, \frac{1}{N} \max_{2 \le m \le N}\, \vert r(1,m)-(m-1) \vert=0 \right)\;=\;1\,,
\end{equation*}
and thus the assumption~\eqref{eq32} holds almost surely. Moreover, \eqref{LLN} also includes 
ergodic resistance sequences
by Birkhoff’s Ergodic Theorem (cf.\ \cite[Theorem 7.2.1]{Durrett}). 
 In this paper, except Section \ref{sec:proofgeneraleig} we always assume \eqref{LLN} when without explicit statement.

\end{rem}

	\begin{theorem}\label{th:gapshapeder} 
		If the condition \eqref{LLN} on the resistances holds,
		we have 
		\begin{equation}\label{gap:asym}
		\lim_{N\to \infty} \frac{N^2 \gap_N }{\pi^2}\;=\;1 \, .
		\end{equation}
		Furthermore, concerning the shape and (weighted) derivative of the eigenfunction $g_1$ with $g_1(1)\colonequals 1$  corresponding to the spectral gap, i.e.\ $\Delta^{(c)} g_1=-\gap_N \cdot g_1$ and  setting 

\begin{equation*}
h(x) \colonequals \cos\left(\frac{\pi(x-1/2)}{N}\right), \;\;\; \forall \, x \in \lint 1, N \rint \, ,
\end{equation*}				
	we have		
\begin{align}
		\lim_{N\to \infty}\; \sup_{x\in \lint 1,\, N\rint}\; &\left| g_1(x) - h(x) \right|\;=\;0\,, \label{eigfun:shape}\\	
\lim_{N \to \infty}\; \sup_{x \in \lint 1,\, N \rint}\; &\left \vert  N (c \nabla g_1)(x)- N(\nabla h)(x) \right\vert \; = \; 0 \,.	\label{approx:der2eigfuns}
\end{align}

	\end{theorem}

\begin{rem}\label{rem:eigshape}
Given $(c(x, x+1))_{1\le x<N}$, 
 the condition $g_1(1)=1$ together with $\gl_1$ already determines the function $g_1$, see  \eqref{grad:g}.
Theorem \ref{th:gapshapeder}--\eqref{approx:der2eigfuns} is a key input for sharp bounds on the mixing time for simple exclusion process, see \cite[Theorem 2.6]{yang2024cutoffSEP}.
 Moreover,
the methods 
in Theorem  \ref{th:gapshapeder} also work for the other $j$-monotone eigenfunctions under the assumption \eqref{LLN}. Precisely,  with $K_0 \in \bbN$ being any prefixed constant,  for all $1 \le j \le K_0$, $g_j(1) \colonequals 1$ and $h_j(x) \colonequals \cos\left(j \pi (x-1/2)/N \right)$, we have
 
\begin{equation}\label{shape:eigfunRC}
\begin{aligned}
&\lim_{N \to \infty} \vert \gl_j N^2/\pi^2-j^2\vert\;=\;0\,,\\
&	\lim_{N\to \infty}\; \sup_{x\in \lint 1,\, N\rint}\; \left| g_j(x)  - h_j(x) \right|\;=\;0\,,\\
&\lim_{N\to \infty}\; \sup_{x\in \lint 1,\, N\rint}\; N\left| (c\nabla g_j)(x)  - (\nabla h_j)(x) \right|\;=\;0\,.\\	
	\end{aligned}
\end{equation} 

\end{rem}

\begin{rem}
If  $(c(x, x+1))_x$ are IID with the law of $1/c(x, x+1)$ being in the $\ga-$attraction domain where $\ga \in (0, 1)$, then $\Gap_{N}$ is of order $N^{-1-1/\ga}$ up to a slowly  varying function in $N$, whose rescaling limit distribution is non-degenerated  \cite{Faggionato2012spectral}. 
\end{rem}

\subsection{Intuition for Theorem \ref{th:gapshapeder}}
 Setting  (quite arbitrarily) $c(N,N+1)=1$, for $\gl>0$,  
let $f^{\gl}: \lint 0, N+1\rint \mapsto \bbR$ be defined by $f^{\gl}(0)=f^{\gl}(1)=1$
and for $x\in \lint 1,N\rint$,
\begin{equation}\label{recurf}
f^{\gl}(x+1)\;=\;f^{\gl}(x)+ \frac{1}{c(x,x+1)}\left[ (c\nabla f^{\gl})(x) -\gl f^{\gl}(x)\right]\,.
\end{equation}
Note  that (the restriction to $\lint 1,N\rint$ of) $f^{\gl}$ is an eigenfunction of $\Delta^{(c)}$ if and only if 
\begin{equation}\label{bdcd:f}
f^{\gl}(N+1)\;=\;f^{\gl}(N)\,,
\end{equation}
since $f^{\gl}$ satisfies $(\Delta^{(c)} f)(x)=-\gl f^{\gl}(x)$ for $x<N$ by construction. 
Up to a multiplicative factor, all eigenfunctions are of this type,  since  there is no eigenfunction satisfying $f^{\gl}(1)=0$  or $f^{\gl}(N)=0$ (explained below \eqref{grad:g}).
 For $\gl>0$ and $x\in \lint 1,N+1\rint$, we set 
\begin{equation}\label{def:blambda}
b(\gl,x)\; \colonequals \;- \frac{(c\nabla f^{\gl})(x)}{f^{\gl}(x-1)}
\end{equation}
with the convention that  $b(\gl,x)=\overline \infty$ if $f^{\gl}(x-1)=0$, 
and consider $\overline \bbR= \bbR \cup \{\overline \infty\}$ to be the Alexandrov compactification of $\bbR$.
By \eqref{recurf} and \eqref{def:blambda}  we deduce that
\begin{equation}\label{recur:b}
b(\gl,x+1)\;=\;\frac{b(\gl,x)}{1-c(x-1,x)^{-1}b(\gl,x)} +\gl\,.
\end{equation}
Setting $B^{(N)}(x) \colonequals b(\gl,x) N$ and $\gl \colonequals \ga/N^2$, by \eqref{recur:b} we have
the recursion
\begin{equation}\label{rel:scale}
B^{(N)}\left(x+1\right)\;=\;\frac{B^{(N)}\left(x\right)}{1- N^{-1} r^{(N)}(x-1,x)B^{(N)} (x)}+\frac{\alpha}{N}\,,
\end{equation}
which starts from $B^{(N)}(1) \colonequals 0$. Note that  by \eqref{bdcd:f} and \eqref{def:blambda},  $\gl$ is an eigenvalue if and only if $B^{(N)}(N+1)=0$. 
Through the paper, we drop the dependence of $\gl$ and the superscript  ``$(N)$'' in  $B^{(N)}(x)$ for ease of notations.
We now provide an intuition for Theorem \ref{th:gapshapeder}--\eqref{gap:asym}: note that \eqref{rel:scale} is equivalent to the following
\begin{equation}\label{resscale:B}
N \left[ B\left( x+1 \right)-B\left( x \right) \right]
\;=\;  \frac{r(x-1, x)B\left( x \right)^2 }{1-B\left( x\right) r(x-1,x) N^{-1}}+\ga\,,
\end{equation}
whose asymptotic ODE can be described by
\begin{equation*}
\begin{cases}
\frac{\dd y}{\dd x}\;=\;  y^2+\ga, \quad x \in [0, \, 1]\,,\\
y(0)\;=\;0\,,
\end{cases}
\end{equation*}
 with its unique solution as 
\begin{equation*}
y(x)\;=\;\sqrt{\ga} \tan \left( \sqrt{\ga} \cdot x \right)\,.
\end{equation*}
Furthermore, in Section \ref{sec:gapestimate} we know that the sequence $(B(x))_{x \le N+1}$ corresponding to the $i$th eigenvalue is well approximated by the first $i$ branches of the tangent function. Concerning the spectral gap, 
since we require $y(1)=0$ and no any other zero in the interval $(0, 1)$,    we obtain  $\ga=\pi^2$. A similar reasoning provides $\gl_j=(1+o(1))j^2 \pi^2/N^2$ in \eqref{shape:eigfunRC}.
In the remaining of the paper,  we turn the intuition into a rigorous proof.

\subsection{Related literatures}

\subsubsection{One-dimensional random conductance models}
 A qualitative spectral framework for random walks in random environments was developed by Faggionato \cite{Faggionato2012spectral}. Using the theory of Krein–Feller operators and generalized second-order differential operators, \cite{Faggionato2012spectral} establishes convergence of eigenvalues and eigenfunctions under ergodicity assumptions on the resistance sequence, primarily for Dirichlet boundary conditions. These results are particularly well suited to the study of subdiffusive trap and barrier models and allow for very general random measures.

However, the conclusions in \cite{Faggionato2012spectral} are largely qualitative: they identify limiting operators and establish convergence, but do not provide sharp asymptotics, explicit eigenfunction profiles, or uniform control of discrete fluxes.
Furthermore, while Faggionato’s theorems focus largely on Dirichlet boundary conditions, the present  paper study utilizes Neumann boundary conditions. This shift is theoretically supported by the Dirichlet-Neumann bracketing technique  \cite[Section 2.2]{Faggionato2012spectral},  which allows for the comparison of spectra across different constraints and justifies the transition to the homogenized limit in this specific context.

By contrast, the present work provides explicit asymptotic formulas for both the spectral gap and the associated principal eigenfunction under a minimal asymptotic homogeneity assumption on the resistances. 
In particular, we show that the eigenfunction converges uniformly to the cosine profile and that its weighted discrete derivative converges as well. These quantitative results are obtained via an elementary perturbative approach and do not rely on operator-theoretic machinery.

\subsubsection{Higher dimensions and homogenization}
 In dimensions $d \ge 2$, Boivin and Depauw \cite{Boivin2003Spectral} established homogenization of the spectrum under uniform ellipticity, which was later extended to weaker assumptions in \cite{Neukamm2017homogenization, Flegel2019Homogenization}. In regimes with heavy-tailed conductances, anomalous spectral behavior and different scaling limits were identified, see  \cite{Flegel2018heavy}. These works focus primarily on convergence of eigenvalues and eigenfunctions at the macroscopic level. 
  Some other properties of the Random Conductance Model are also investigated,  like quenched invariance principal \cite{Andres2015Invariance} on $\bbZ^d$, the speed \cite{Berger2013speed, Berger2019speed} on $\bbZ^2$, heat kernel estimates  \cite{Barlow2004supercritical,  Mathieu2004Isoperimetry} on percolation clusters and so on. We refer to \cite[Section 1.3]{Dario2021quantitative} for comments on recent progress on the random conductance model and the stochastic homogenization method.   When the conductances on each edge are homogenously equal to one, concerning general graphs we refer to \cite{Friedman1993geometric}
for geometric insight into the eigenvectors. Additionally, we point to \cite[Section 1.1.1]{Gross2024eigenvalues}
for related works about the relations between eigenvalues and the heat  kernel.

\subsubsection{Connections to stochastic homogenization}

Our results are also related to the broader theory of stochastic homogenization. Classical works based on 
G-convergence and compactness methods \cite{Jikov1994Homogenization} establish convergence of spectra for elliptic operators with rapidly oscillating coefficients but do not yield convergence rates or explicit profiles.

More recently, quantitative stochastic homogenization has achieved optimal rates for correctors and homogenized coefficients, notably in the works of Gloria, Neukamm, and Otto
\cite{Gloria2011optimal, Gloria2012error, Gloria2015Quantification} and Armstrong, Kuusi, and Mourrat \cite{Armstrong2019homogen}. Spectral fluctuations and convergence rates for eigenvalues were investigated in \cite{Duerinckx2022eigen}. These results concern continuum operators and rely on large-scale regularity techniques that are quite different in nature from the discrete, one-dimensional approach adopted here.

\subsection*{Organization}
Section \ref{sec:proofgeneraleig} is about the proof of Proposition \ref{prop:eigfunc} concerning monotonicities of eigenfunctions in arbitrary disordered setup.
 Section \ref{sec:gapestimate} is devoted to the estimate of the spectral gap under the assumption \eqref{LLN}.
 Section \ref{sec:eigshapeder} is devoted to prove the shape and derivative of the principal eigenfunction under the assumption \eqref{LLN}.

\subsection*{Acknowledgments}
S.Y.\ is very grateful to  Hubert Lacoin for suggesting this problem and  thanks  Tertuliano Franco and Hubert Lacoin for enlightening and insightful discussions. 
  Moreover, he also thanks sincerely Gideon Amir,  Chenlin Gu, Gady Kozma, and Leonardo Rolla for helpful discussions. 
  S.Y.\ is supported by CNPq 401314/2025-1. 
   Partial of this work was done during his stay in Universidade de São Paulo supported by FAPESP 2023/12652-4, and in
 Bar-Ilan University supported by  Israel Science Foundation grants 1327/19
and 957/20.

\section{Proof of Proposition \ref{prop:eigfunc}}\label{sec:proofgeneraleig}

Note that in this section we prove Proposition \ref{prop:eigfunc} only with the assumption $c(x,x+1)>0$ for all $x \in \lint 1, N-1\rint.$

\subsection{Proof for Item (1) of Proposition \ref{prop:eigfunc}}
	Since the uniform probability measure $\mu_{N}$ on $\lint 1, N \rint$
	satisfies the detailed balance condition, we know that (cf.\
	 \cite[Proposition 3.1 in Chapter 1]{Yang2021thesis}) there exists $N$ distinct eigenfunctions $\{ g_i\}_{i=0}^{N-1}$ with corresponding real eigenvalues $\{ -\gl_i\}_{i=0}^{N-1}$ of the generator $\Delta^{(c)}$
	such that $g_0=\ind$ and
	\begin{equation}\label{sets:eigfuns}
	\begin{cases}
	0\;=\;\gl_0\,<\,\gl_1\,\le\, \gl_2\,\le\, \cdots \,\le\, \gl_{N-1}\,,\\
	\Delta^{(c)} g_i\;=\;-\gl_i g_i \text{ and } g_i(1)\;=\;1\,, \quad &\forall\, i \in \lint 0, N-1 \rint\,,\\
	\frac{1}{N}\sum_{x=1}^N g_i(x)g_j(x)= C_{i,j} \delta_{i,j}\,, \quad &\forall\, i,j \in \lint 0, N-1 \rint\,,
	\end{cases}
	\end{equation}
	where $\delta_{i,j}$ is the Kronecker delta, and $(C_{i,i})_{i}$ are some positive constants ($C_{i,j}$ with $i \neq j$ does not play a role).
	By \eqref{randwalk:generator} and $(\Delta^{(c)}g_i)(x)=-\gl_i g_i(x)$ for all $x \in \lint 1, N \rint$,
	we have
	\begin{equation}\label{grad:g}
	(c\nabla g_i)(x+1)\;=\;- \gl_i\sum_{k=1}^x g_i(k)\,.
	\end{equation}
By \eqref{grad:g}, given $(c(x, x+1))_{1\le x<N}$ we know that  $\gl_i$ and  $(g_i(k))_{1\le k \le x}$ determines $g_{i}(x+1)$, and thus   $\gl_i$ determines $g_i$.  Here also explains  $g_i(1) \neq 0$, otherwise $g_i \equiv 0$. Similarly,  $g_i(N)\neq 0$  for all $i \in \lint 0, N-1\rint$.
	Since $g_i \neq g_j$ for $i \neq j$, we obtain that $\gl_i \neq \gl_j$ for all $ i \neq j$.
	Therefore, we have proved 
	\eqref{eigv:distinct}.
	\qed

\subsection{Strict monotonicity of the principal eigenfunction}
	We postpone  the proof of the $j$-monotonicity of $g_j$ to the end of this section. So, for the moment, let us assume that $g_1$ is decreasing. 
	To lighten notations in the sequel,  $g_1$  is simply  denoted by $g$.
	We now prove that $g$ is indeed strictly decreasing, i.e.\
	\begin{equation}\label{delta_minimum}
	\delta_{\min} \;\colonequals\; \min_{2\le x \le N} \left[ g(x-1)-g(x)\right]\,>\,0\,.
	\end{equation}
	Suppose the claim does not hold, and we argue by contradiction. By $g(1)=1$,  we observe that
	 $g(1)>g(2)$ by \eqref{grad:g}, and    $g(N)<0$ by  $\sum_{i=1}^N g(i)=0$ and the fact that $g$ is decreasing.	
Furthermore, $g(N-1)>g(N)$ by $(c \nabla g)(N)=\gl_1 g(N)$. Then there exists a smallest integer  $x_0 \in \lint 2, N-2 \rint$ such that
	\begin{equation*}
	\delta_{\min}\;=\;g(x_0)-g(x_0+1)\;=\;0\,.
	\end{equation*}
	By \eqref{grad:g} and $\sum_{i=1}^N g(i)=0$, we have $\sum_{i=1}^{x_0} g(i)=0$ and $\sum_{i=x_0+1}^N g(i)=0$. Plugging $g$ into \eqref{def:gap}, we obtain
	\begin{equation*}	
	\gl_1  \;=\;  \frac{\sum_{k=1}^{x_0-1} c(k,k+1)\left( g(k+1)-g(k) \right)^2 + \sum_{k=x_0+1}^{N-1} c(k,k+1)\left( g(k+1)-g(k) \right)^2}{\sum_{k=1}^{x_0}g(k)^2+ \sum_{k=x_0+1}^{N}g(k)^2}\,.
	\end{equation*}
	Then we have
	\begin{equation}\label{varform:gap}
	\frac{\sum_{k=1}^{x_0-1} c(k,k+1)\left( g(k+1)-g(k) \right)^2}{\sum_{k=1}^{x_0}g(k)^2} \le \gl_1 \;\text{ or  }\; \frac{\sum_{k=x_0+1}^{N-1} c(k,k+1)\left( g(k+1)-g(k) \right)^2}{\sum_{k=x_0+1}^{N}g(k)^2} \le \gl_1\,.
	\end{equation}
	Without loss of generality, suppose  that the first inequality in \eqref{varform:gap} holds.
We	define
	a function $\tilde g: \lint 1, N \rint \mapsto \bbR$ by
	\begin{equation*}
	\tilde g(x) \colonequals
	\begin{cases}
	g(x) &\text{ if } x \le x_0\,,\\
	g(x_0) &\text{ if } x> x_0\,.
	\end{cases}
	\end{equation*}
	Now	we claim that $g(x_0)<0$. In fact, by the definition of $x_0$ and \eqref{grad:g}, we obtain  that
	\begin{equation*}
	\sum_{k=1}^{x_0-1}g(k)>0\quad\text{ and }\quad \sum_{k=1}^{x_0}g(k)=0\,,
	\end{equation*}
	which implies $g(x_0)<0$. Since $\sum_{k=1}^{x_0}g(k)=0$, evaluating the mean and variance of $\tilde{g}$ with respect to $\mu_N$, we have
	\begin{equation*}
	\mu_N(\tilde g)= \frac{N-x_0}{N} g(x_0) \quad\text{ and }\quad \Var_{\mu_N}(\tilde g)=\frac{1}{N} \sum_{x=1}^{x_0}g(x)^2 +\frac{x_0(N-x_0)}{N^2} g(x_0)^2\,.
	\end{equation*}
	Therefore, plugging $\tilde g$ in the variational formula \eqref{def:gap} and comparing with the first inequality of 
	 \eqref{varform:gap}, we obtain 
	\begin{equation*}
	\frac{-\langle \Delta^{(c)}\tilde g, \tilde g \rangle_{\mu}}{\Var_{\mu}(\tilde g)}\;=\;\frac{\frac{1}{N} \sum_{x=1}^{x_0-1}c(x, x+1)\left(g(x+1)-g(x) \right)^2}{\frac{1}{N} \sum_{x=1}^{x_0}g(x)^2 +\frac{x_0(N-x_0)}{N^2} g(x_0)^2}\;<\; \gl_1
	\end{equation*}
	which is a contradiction to the variational formula \eqref{def:gap} for the spectral gap. The other case is treated similarly. Therefore, we have
	$ \delta_{\min}>0.$
	\qed

\subsection{Proof for the $j$-monotonicity property in Proposition \ref{prop:eigfunc}}

Concerning	the monotonicities of the $j$th-eigenfunction stated in Proposition~\ref{prop:eigfunc}, 
	the proof in \cite{Miclo2008eigenfunctions} is in an analytical viewpoint. Here we provide an alternative self-contained proof  based on a more geometric viewpoint, which will be an essential input for a sharp estimate on the spectral gap.

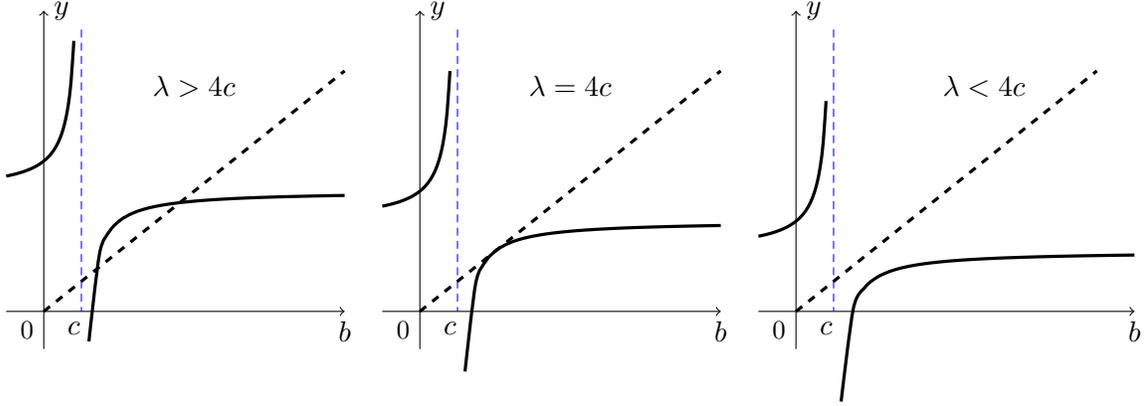
\begin{figure}[H]
	\centering
	\begin{tikzpicture}[scale=0.5,smooth];
	\draw[->] (0,-1)--(0,8) node[right]{$y$};
	\node at (4,6){$\gl > 4c$};
	\draw[->] (-1,0)--(8,0) node[anchor=north]{$b$};

	\draw[black = solid,  very thick ] plot
	[domain=-1:0.8](\x,{(\x/(1-\x)+5)*0.8});

	\draw[black = solid,  very thick ] plot
	[domain=1.2:8](\x,{(\x/(1-\x)+5)*0.8});;

	\draw[dashed, black = solid,  very thick ] plot
	[domain=0:8](\x,{\x*0.8});

	\draw (0,0) node[anchor=north east]{\small $0$};

	\draw[densely dashed,blue] (1,7.5) -- (1,0);
	\node at (0.8,0)[below]{$c$};

	\draw[->] (10,-1)--(10,8) node[right]{$y$};
	\node at (14,6){$\gl = 4c$};
	\draw[->] (9,0)--(18,0) node[anchor=north]{$b$};

	\draw[black = solid,  very thick ] plot
	[domain=9:10.8](\x,{((\x-10)/(11-\x)+4)*0.8});

	\draw[black = solid,  very thick ] plot
	[domain=11.2:18](\x,{((\x-10)/(11-\x)+4)*0.8});;

	\draw[dashed, black = solid,  very thick ] plot
	[domain=10:18](\x,{(\x-10)*0.8});

	\draw (10,0) node[anchor=north east]{\small $0$};

	\draw[densely dashed,blue] (11,7.5) -- (11,0); \node at (10.8,0)[below]{$c$};

	\draw[->] (20,-1)--(20,8) node[right]{$y$};
	\node at (25,6){$\gl < 4c$};
	\draw[->] (19,0)--(29,0) node[anchor=north]{$b$};
	\draw[black = solid,  very thick ] plot
	[domain=19:20.8](\x,{((\x-20)/(21-\x)+3)*0.8});

	\draw[black = solid,  very thick ] plot
	[domain=21.2:29](\x,{((\x-20)/(21-\x)+3)*0.8});

	\draw[dashed, black = solid,  very thick ] plot
	[domain=20:28](\x,{(\x-20)*0.8});

	\draw (20,0) node[anchor=north east]{\small $0$};

	\draw[densely dashed,blue] (21,7.5) -- (21,0); \node at (20.8,0)[below]{$c$};

	\end{tikzpicture}
	\caption{In the figures above,  solid lines depict the function $\Xi^{(c)}(b,\gl)$ with $\gl>0$ fixed,  the black dashed lines stand for $y=b$, and the blue dashes lines are $b=c$. There are  three cases about the number of fixed points of the map $b \mapsto \Xi^{(c)}(\gl, b)$  according to the relation between $\gl$ and $4c$ as shown above.}
	\label{fig:fixpoint}
\end{figure}

 From \eqref{recurf} we remark that $c(N,N+1)$ does not play a role for $b(\gl,N+1)$.
Inspired by \eqref{recur:b}, given a fixed $c>0$, we define $\Xi^{(c)}: \bbR\times \overline \bbR \to \overline \bbR$ as
\begin{equation}\label{def:Xi}
\Xi^{(c)}(\gl,b)\;=\;\frac{b}{1-c^{-1}b}+\gl\,.
\end{equation}
The function $b\mapsto \Xi^{(c)}(\gl,b)$ may have zero, one or two fixed points depending on the values of $\lambda$ and $c$, see Figure~\ref{fig:fixpoint}.
If $b\mapsto \Xi^{(c)}(\gl,b)$ has fixed points $b_1$ and $b_2$ (not necessarily distinct) such that  $ b_1\le b_2$, 
we define $I(\gl,c)=[b_1,b_2]$.  Otherwise, set $I(\gl,c)=\varnothing$ .
We then define the ``angle mapping'' function
\begin{equation}\label{def:varphi}
\varphi(c,\gl,\theta)\;:=\; \inf\{ \theta'\ge \theta
+ \pi \ind_{I(\gl,c)}(\tan \theta)
\ : \ \tan \theta'=  \Xi^{(c)}(\gl,\tan \theta) \}\,,
\end{equation}
with the convention that $\tan(\pi/2+k\pi)=\overline \infty$ for $k\in \bbZ$.
We now recursively define  an ``angle'' $\theta(\gl,x)$ by setting $\theta(\gl,1)=0$ and for $x \in  \lint 1, N \rint$,
\begin{equation}\label{def:maptheta}
\theta(\gl,x+1)\;:=\;  \varphi(c(x-1,x), \gl,\theta(\gl,x))\,.
\end{equation}
For $\gl=0$, we define $\theta(0,x)=0$ for all $x \in \bbN$, and  the map $b \mapsto \Xi^{(c)}(0,b)$  has only one fixed point $b=0$. We are going to prove the following properties for the function $\gl\mapsto \theta(\gl,x)$.
\begin{lemma}\label{lema:anglecont}
	For fixed $c, \gl>0$, the map   $\theta \mapsto \varphi(c,\gl, \theta)$ in \eqref{def:varphi} is continuous and strictly increasing.
\end{lemma}

\begin{proof}
	For $b\le 0$, we have $\Xi^{(c)}(\gl,b)= \frac{b}{1-c^{-1}b}+\gl>b$. Therefore, if there exists fixed points $b_1,b_2$ (not necessarily distinct) for the map $b \mapsto \Xi^{(c)}(\gl,b)$, we must have $0<b_1 \le b_2$.
	Referring to Figure \ref{fig:fixpoint} for a graphic description of the function $\Xi^{(c)}$ (with $\gl$ and $ c$ fixed), we observe that
	\begin{itemize}
		\item If there exists at least one fixed point of the map $b \mapsto \Xi^{(c)}(\gl,b)$, then we have $\gl \ge 4c$. Moreover, if $\gl>4c$, there are two distinct fixed points. While if $\gl=4c$, there is only one fixed point.
	\end{itemize}
Moreover, note that in \eqref{def:varphi} when $\tan \theta=c$, $\tan \theta' \in \{ \pm \infty\}$ gives the same output in $\varphi$.
	Therefore we can identify $\{ \pm \infty  \}$ to be $\overline \infty$, and then $\Xi^{(c)}(\gl,\cdot):  \overline{\bbR} \mapsto \overline{\bbR}$ is continuous.

	Let $\theta_1$ be an arbitrary point. We first deal with the continuity of the map $\theta \mapsto \varphi(c, \gl, \theta)$.
	If $\theta_1$ satisfies $\tan \theta_1 \neq b_1,b_2$,
	we can choose $\delta>0$ sufficiently small such that if $\vert \theta- \theta_1 \vert < \delta$, we have
	$\ind_{I(\gl,c)}(\tan \theta)= \ind_{I(\gl,c)}(\tan \theta_1).$
	That $\varphi(c, \gl, \cdot)$ is  continuous at $\theta_1$ follows from the continuity of the maps $\Xi^{c}(\gl, \cdot)$ and  $\arctan(\cdot)$.

	Now we deal with the case $\theta_1$ satisfying $\tan \theta_1 = b_1$ or $b_2$. Without lost of generality, we assume $\tan \theta_1=b_1$. If $b_1<b_2$, we have to analyze it according to whether $\theta$ approaches $\theta_1$ from the left or from the right.  If   $0 \le \theta-\theta_1 < \delta$ is such that  $\ind_{I(\gl,c)}(\tan \theta)=\ind_{I(\gl,c)}(\tan \theta_1)$, then we can conclude that  $\varphi(c,\gl,\cdot)$ is right-continuous at $\theta_1$.
	Now let $0 <\theta_1-\theta< \delta$, and then we have $\Xi^{(c)}(\gl, \tan \theta)< \tan \theta$. For any given $\epsilon>0$, we can choose $\delta=\delta(\gep)$ such that $  \varphi(c, \gl, \theta_1)-\gep<\varphi(c, \gl, \theta) \le \varphi(c, \gl, \theta_1) $ due to the continuity of $\tan(\cdot)$ function and that $\varphi(c, \gl, \tan \theta)$ is in next branch of the tangent  function.
	Thus $\varphi(c, \gl,\cdot)$ is left continuous at $\theta_1$.

	For the case $\theta_1$ satisfying $\tan \theta_1=b_1=b_2$, we can argue similarly. Therefore, we prove the continuity of the map $\theta \mapsto \varphi(c,\gl, \theta)$.

\medskip

	Now we move to prove that this map is also strictly increasing in $\theta$. It is sufficient to show that for any arbitrary $\theta_1$, there exists  $\delta=\delta(\theta_1)$ such that the map $\varphi(c, \gl, \cdot): (\theta_1-\delta, \theta_1+\delta) \mapsto \bbR $ is strictly increasing.
	We first deal with the case $\tan \theta_1=c$. As $c$ is not a fixed point of the map $b \mapsto \Xi^{(c)}(\gl,b)$, we can choose $\delta=\delta(\theta_1)$ sufficiently small such that the map $b \mapsto \Xi^{(c)}(\gl,b)$ has no fixed point in the domain
	$  (\tan(\theta_1-\delta), c) \cup (c, \tan( \theta_1+\delta)).$
	Then we have for all $\theta \in (\theta_1-\delta, \theta_1+\delta)$,
	\begin{equation*}
	\ind_{I(\gl,c)}(\tan \theta)\;=\;\ind_{I(\gl,c)}(\tan \theta_1)\;=\;0\,.
	\end{equation*}
	Since $\varphi(c, \gl, \theta_1-)=\varphi(c, \gl, \theta_1+)$, we can use that $b \mapsto \Xi^{(c)}(\gl,b)$ is strictly increasing in $(\tan(\theta_1-\delta), c) \cup (c, \tan( \theta_1+\delta))$ to obtain that $\varphi$ is strictly increasing in $(\theta_1-\delta, \theta_1+\delta)$.

	From now on we deal with the case $\tan \theta_1 \neq c$. If $\tan \theta_1 \neq b_1,b_2$, we can choose $\delta=\delta(\theta_1)$ sufficiently small such that for all $\theta \in (\theta_1-\delta, \theta_1+\delta)$,
	\begin{equation*}
	\ind_{I(\gl,c)}(\tan \theta_1)\;=\;\ind_{I(\gl,c)}(\tan \theta)
	\end{equation*}
	and the map $b \mapsto \Xi^{(c)}(\gl, b)$ restricted to $(\tan(\theta_1-\delta), \tan(\theta_1+\delta))$ is strictly increasing. Therefore, $\varphi(c, \gl, \cdot)$ is strictly increasing in $(\theta_1-\delta, \theta_1+\delta)$.

	Now we deal with the case $\tan \theta_1=b_1$ or $b_2$ where $b_1 \neq b_2$. Suppose w.l.o.g.\ that $\tan \theta_1=b_1$.    We can choose $\delta=\delta(\theta_1)>0$ sufficiently small such that for all $\theta \in [\theta_1, \theta_1+\delta)$, we have $\ind_{I(\gl,c)}(\tan \theta_1)=\ind_{I(\gl,c)}(\tan \theta)$ and that
	the map $b \mapsto \Xi^{(c)}(\gl, b)$ in $[\tan \theta_1, \tan(\theta_1+\delta))$ is strictly increasing. Therefore, $\varphi(c, \gl, \cdot)$ is strictly increasing in $[\theta_1, \theta_1+\delta)$.
	While for $\theta \in (\theta_1-\delta, \theta_1]$, we know that
	the map $b \mapsto \Xi^{(c)}(\gl, b)$ restricted to $(\tan (\theta_1-\delta), \tan \theta_1)$ is strictly increasing, and also that  $\ind_{I(\gl,c)}(\tan \theta_1)=1$ and $\ind_{I(\gl,c)}(\tan \theta)=0$ for $\theta \in (\theta_1-\delta, \theta_1)$. Therefore,
	$\varphi(c, \gl, \cdot)$ is strictly increasing in $(\theta_1-\delta, \theta_1]$.

	About the case $\tan \theta_1=b_1=b_2$, for $\delta=\delta(\theta_1)>0$ sufficiently small, the map $\varphi(c, \gl, \cdot): (\theta_1-\delta, \theta_1]\mapsto \bbR $ is clearly strictly increasing with the same reasoning as above. For $\theta \in [\theta_1, \theta_1+\delta)$, the map $\Xi^{(c)}(\gl, \cdot): [\tan \theta_1, \tan(\theta_1+\delta)) \to \bbR $ is strictly increasing, and $  \Xi^{(c)}(\gl, \tan \theta)< \tan \theta$ for $\theta \in (\theta_1, \theta_1+\delta)$. Therefore, $\varphi (c, \gl, \theta)$ is in next branch due to the $\tan$ function, and hence $\varphi (c, \gl, \cdot)$ is  strictly increasing in $[\theta_1, \theta_1+\delta)$.
	Thus we have proved that  $\varphi(c, \gl, \cdot)$ is strictly increasing in $\theta$.
\end{proof}

\begin{lemma}\label{lema:lambdamonocont}
	For fixed $c, \theta>0$, the map $ \gl \mapsto \varphi(c, \gl, \theta)$  is strictly increasing and  uniformly continuous  in $\theta$.
\end{lemma}

\begin{proof}
	Note that $\Xi^{(c)}(b, \cdot)$ is  strictly  increasing in $\gl$.  We  claim that $I(\gl_1, c)  \subseteq I(\gl_2, c)$ for $\gl_1< \gl_2$, which is enough to conclude that $\varphi(c,\cdot, \theta)$ is strictly  increasing in $\gl$. It is sufficient to deal with $\gl_1 \ge 4c$, because otherwise there is no fixed points (see Figure~\ref{fig:fixpoint}) and the claim holds trivially. By \eqref{def:Xi}, we have
	\begin{equation}\label{fixdomain}
	I(\gl,c)\;=\;\left[\frac{\gl-\sqrt{\gl^2-4\gl c}}{2},\, \frac{\gl+\sqrt{\gl^2-4\gl c}}{2}  \right] \subset (0, \infty)\,.
	\end{equation}
	Furthermore, for $\gl \ge 4c,$ the function $\gl \mapsto \frac{\gl-\sqrt{\gl^2-4\gl c}}{2}$  is strictly decreasing, and the function $\gl \mapsto \frac{\gl+\sqrt{\gl^2-4\gl c}}{2} $ is strictly increasing. Therefore, $\varphi$ is strictly  increasing in $\gl$.

	We now move to prove that $\varphi$ is continuous in $\gl$ uniformly in $\theta$.  Fix $\gl_1 \ge 0$. By \eqref{fixdomain}, if $\tan \theta \neq b_1(\gl_1,c), b_2(\gl_1,c)$, there exists $\delta>0$ such that if $ \vert \gl_1-\gl \vert< \delta$, we have $\ind_{I(\gl,c)}(\tan \theta)=\ind_{I(\gl_1,c)}(\tan \theta)$.
	Since both $ \gl \mapsto \Xi^{(c)}(b, \gl)$ and $\arctan: \bbR \to (-\frac{\pi}{2}+n \pi, \frac{\pi}{2}+n \pi)$ where $n \in \bbN$ are $1$-Lipschitz continuous,  then $\varphi$ is continuous at $\gl_1$ uniformly in $\theta$.

	Now we assume w.l.o.g.\ that $\tan \theta =b_1(\gl_1, c)<b_2(\gl_1, c)$. If $0\le \gl-\gl_1 \le \delta$, we have $\ind_{I(\gl,c)}(\tan \theta)=\ind_{I(\gl_1,c)}(\tan \theta)=1$ and then $\varphi$ is right-continuous at $\gl_1$ uniformly in $\theta$ by the $1$-Lipschitz continuity of the maps $ \gl \mapsto \Xi^{(c)}(b, \gl)$ and $\arctan: \bbR \to (-\frac{\pi}{2}+n \pi, \frac{\pi}{2}+n \pi)$.
	On the other hand, if $0< \gl_1-\gl< \delta$, we observe that
	\begin{align*}
	&\Xi^{(c)}(\gl_1, \tan \theta)-\Xi^{(c)}(\gl, \tan \theta)\;=\;\gl_1-\gl\quad \text{ and }\\
	&0 \;\le\; \arctan \Xi^{(c)}(\gl_1, \tan \theta)-\arctan \Xi^{(c)}(\gl, \tan \theta) \;\le\; \gl_1-\gl
	\end{align*}
	where we pick $\arctan: \bbR \mapsto (-\frac{\pi}{2}+n \pi, \frac{\pi}{2}+n \pi]$ with $n \in \bbN$ chosen such that $\varphi (c, \gl_1, \theta) \in (-\frac{\pi}{2}+n \pi, \frac{\pi}{2}+n \pi]$.
	Moreover, we have $\Xi^{(c)}(\gl, \tan \theta)<\tan \theta$, see the leftmost picture in Figure~\ref{fig:fixpoint}.  By the definition of $\varphi$ in \eqref{def:varphi}, $\varphi(c, \gl, \theta)$ lies in next branch of the tangent function.
	Besides, since $\tan \theta$ is a fixed point (for the map $\Xi^{(c)}(\cdot,\lambda_1)$), we have that $\ind_{I(\gl_1,c)}(\tan \theta)=1$ and then $\varphi(c, \gl, \theta)$ and $\varphi(c, \gl_1, \theta)$ are in the same branch of the tangent function. Therefore,
	$ 0 \le \varphi(c, \gl_1, \theta)-  \varphi(c, \gl, \theta) \le \gl_1-\gl, $
	i.e.\ $\varphi(c, \cdot, \theta)$ is left-continuous at $\gl_1$ uniformly in $\theta$.

	For the case $\tan \theta=b_1=b_2$, the proof is similar.
	As $\gl_1$ is arbitrary, we conclude the proof.

\end{proof}

\begin{lemma} \label{lema:jointcont}
	For fixed $c>0$, the map $(\gl, \theta) \mapsto \varphi(c, \gl, \theta)$ is jointly continuous.
\end{lemma}

\begin{proof}
	Fixing $(\gl_1, \theta_1)$, we have
	\begin{equation*}
	\varphi(c, \gl, \theta)-\varphi(c, \gl_1, \theta_1)\;=\;\varphi(c, \gl, \theta)-\varphi(c, \gl_1, \theta)+\varphi(c, \gl_1, \theta)-\varphi(c, \gl_1, \theta_1).
	\end{equation*}
	Then we conclude the proof by Lemma \ref{lema:lambdamonocont} and Lemma \ref{lema:anglecont}.
\end{proof}

\begin{lemma}\label{lema:thetacont}
For $\gl > 0$,	the map
	$\gl \mapsto \theta(\gl,x)$, defined in \eqref{def:maptheta}, is continuous and strictly increasing for any $x \in \lint 2, N \rint$.
\end{lemma}
\begin{proof}
	We first deal with the continuity and monotonicity of $\gl \mapsto \theta(\gl, 2)$. We recall that $\theta(\gl,1)=0$ for all $\gl \ge 0$. The fact  that the map $\gl \mapsto \theta(\gl,2)$ is strictly increasing and continuous  follows  from that 
	 the map $\gl \mapsto \varphi(c(0,1), \gl, 0)=\arctan(\gl)\in (0, \pi/2)$ is strictly increasing and continuous,  where $c(0,1)>0$ is arbitrary.

	Now we proceed by induction,  assuming that $\gl \mapsto \theta(\gl,x)$ is continuous and strictly increasing. Since
	$\theta(\gl, x+1)= \varphi(c(x-1,x), \gl, \theta(\gl,x))$,
	by Lemma \ref{lema:jointcont}
	we obtain the continuity of the map $\gl \mapsto \theta(\gl,x+1)$.
	By Lemma \ref{lema:anglecont}
	and Lemma \ref{lema:lambdamonocont}, we know that $\theta(\gl, x+1)= \varphi(c(x-1,x), \gl, \theta(\gl,x))$ is strictly increasing in $\gl$. Therefore, we conclude the proof.
\end{proof}

Now note that  $f^{\gl}$ is an eigenfunction if and only if 
$\theta(\gl,N+1)$ is a multiple of $\pi$.
Lemma \ref{lema:thetacont}, together with the fact that  $\theta(0,N+1)=0$, implies
\begin{equation*}
f^{\gl} \text{ is an eigenfunction } \quad \Leftrightarrow\quad  \theta(\gl,N+1)=k\pi \; \text{ for some  }  k\in \lint 0,N-1\rint\,.
\end{equation*}
Let $\gl_k>0$ denote the unique number satisfying $\theta(\gl_k,N+1)=k\pi$, whose existence is assured by Lemma~\ref{lema:thetacont}, and set $f_k:=f^{\gl_k}$.
Let $x_i \in \lint 1, N\rint$  such that $\theta(\gl_k,x_i)\le  i\pi <\theta(\gl_k,x_i+1)$ for $i \in \lint 1, k-1\rint$.

\begin{lemma}\label{lema:numx}
 For $\gl_k$ mentioned above and the associated sequence $(x_i)_i$,  we have $1<x_i < N$, and  $\#\left\{ (x_i)_i\right\}=k-1$.
\end{lemma}

\begin{proof}
Since $\theta(\gl_k, 1)=0$ by definition and $0<\theta(\gl_k, 2)< \frac{\pi}{2}$, we have $x_1 \ge 2$.
Now
suppose there exists one  $x_\ell$ satisfying $x_\ell=N$, then we have
\begin{equation*}
\theta(\gl_k, x_\ell) \;\le\; (k-1)\pi \;<\; \theta(\gl_k, x_\ell+1)\,.
\end{equation*}
  We argue in three cases to obtain a contradiction to $x_\ell=N$.
If $\tan \theta(\gl_k, x_\ell)>0$, then
$(k-2)\pi< \theta(\gl_k, x_\ell) < (k-2)\pi+\frac{\pi}{2}$. As $x_\ell+1=N+1$ and $\theta(\gl_k, N+1)=k \pi$, then $\theta(\gl_k, x_\ell)$ and $\theta(\gl_k, x_\ell+1)$ are in two nonadjacent branches of the $\tan$ function which cannot occur by \eqref{def:varphi}.

If $\tan \theta(\gl_k, x_\ell)=0$, we know that $f_k(N)\neq 0$ (explained below \eqref{grad:g}) and $f_k(N-1)=f_k(N)=f_k(N+1)$, which is a contradiction to $(\Delta^{(c)}f_k)(N)=-\gl_k f_k(N)$.

If $\tan \theta(\gl_k, x_\ell)<0$, since the two fixed points (not necessary distinct) of the map $b \mapsto \Xi^{(c)}(\gl, b)$ are strictly positive (cf. Figure \ref{fig:fixpoint}), we have 
 $\theta(\gl_k, x_\ell+1) \le (k-1)\pi+ \frac{\pi}{2}$. This is a contradiction to $x_\ell+1=N+1$ and $\theta(\gl_k, N+1)=k \pi$.

\medskip

Now we move to show $\# \{(x_i)_i \}=k-1$. By definition, we have $\# \{(x_i)_i \} \le k-1$. Note that
by \eqref{def:Xi} and \eqref{def:varphi}, any two adjacent angles 
 $\theta(\gl_k, y)$ and $\theta(\gl_k, y+1)$ are either in the same branch or adjacent branches of the tangent function. If $\#\{ (x_i)_i\}<k-1$, there exists an adjacent pair of angles $\theta(\gl_k, y)$ and $\theta(\gl_k, y+1)$ in two adjacent branches of the tangent function satisfying for some $i_0 \in \lint 1, k-1\rint$,
 \begin{equation}
i_0 \pi- \frac{\pi}{2} \;\le\; \theta(\gl_k, y) \;\le\; i_0 \pi\,, \; \text{ and }\; (i_0+1) \pi \;<\; \theta(\gl_k, y) \;\le\;  (i_0+1) \pi+ \frac{\pi}{2}\,.
 \end{equation}
 Thus $\tan \theta(\gl_k, y)<0$ and $\theta(\gl_k, y+1)-\theta(\gl_k, y)> \pi$, which is impossible due to any fixed points $b_1(\gl_k, c)>0,$  $b_2(\gl_k, c)>0$ (cf. Figure \ref{fig:fixpoint}).

\end{proof}

Recalling \eqref{def:blambda} and \eqref{def:varphi}, we have
\begin{equation*}
\tan \left( \theta(\gl_k,x) \right)=-\frac{c(x-1,x) \left[f_k(x)-f_k(x-1) \right]}{f_k(x-1)}\,.
\end{equation*}

\begin{lemma}\label{lema:extrema}
	The points $(x_i)^{k-1}_{i=1}$ (or the pair $\{x_i-1,x_i \}$ when $\theta(\gl_k,x_i)=  i\pi$), associated with $\gl_k$, are the local extrema of $f_k$.
\end{lemma}

\begin{proof}
	We first deal with the case $\tan \theta(\gl_k, x_i)=0$, and then $f_k(x_i)=f_k(x_i-1)$. Note that $\tan \theta(\gl_k, 2)>0$, and thus $x_i \ge 3$.
	Since  for all $x \in \lint 2, N-1 \rint$ and $\gl_k >0$,
	\begin{align*}
	(\Delta^{(c)} f_k)(x)& \;=\;c(x-1,x)\left( f_k(x-1)-f_k(x) \right)+c(x,x+1) \left( f_k(x+1)-f_k(x) \right)\\
	&\;=\;-\gl_k f_k(x)\,,
	\end{align*}
	then there is no $x \in \lint 2, N-1 \rint$ satisfying any of the following cases
	\begin{equation}\label{eq2extr}
	\begin{gathered}
	f_k(x-1)\;=\;f_k(x)\;=\;0\,,\\
	0\; <\; f_k(x) \;\le\; \min \left( f_k(x-1), f_k(x+1) \right)\,,\\
	0 \;>\; f_k(x)  \;\ge \;\max \left( f_k(x-1), f_k(x+1) \right)\,.\\
	\end{gathered}
	\end{equation}
If either of the following equality holds:
\begin{equation*}
\begin{gathered}
f_k(x_i-2)\;<\; f_k(x_i-1)\;=\;f_k(x_i)\;<\; f_k(x_i+1)\,,\\
f_k(x_i-2)\;>\; f_k(x_i-1)\;=\;f_k(x_i)\;>\; f_k(x_i+1)\,,
\end{gathered}
\end{equation*}
by analyzing whether $f_k(x_i-1)=f_k(x_i)>0$ or $f_k(x_i-1)=f_k(x_i)<0$, 
 we reach a contradiction to \eqref{eq2extr}, and thus  both $x_i$ and $x_i-1$ are local extrema.

	Now we treat the situation $\tan \theta(\gl_k,x_i)= \overline \infty$, where $f_k(x_i-1)=0$, $\theta(\gl_k,x_i)= i\pi-\frac{\pi}{2}$, and $\ind_{I(\gl_k, c(x_i-1,x_i))}(\tan \theta(\gl_k,x_i) )=0$.
	By the definition of $\varphi$, for all $\lambda > 0, \; c>0,$ and  $ \theta\in \bbR$, if $\ind_{I(\gl,c)}(\tan \theta)=0$ we have 
	\begin{equation}\label{smaller_pi}
	  \varphi( c, \lambda, \theta)-\theta \;\le\; \pi\,. 
	\end{equation}
	The combination of \eqref{smaller_pi} and the assumption $\theta(\gl_k,x_i)\le  i\pi <\theta(\gl_k,x_i+1)$ implies that
	\begin{equation}\label{localextremumden}
\tan \theta(\gl_k, x_i+1)\;=\;	-\frac{c(x_i,x_i+1) \left(f_k(x_i+1)-f_k(x_i) \right)}{f_k(x_i)}\;>\;0\,.
	\end{equation}
	Using $f_k(x_i-1)=0$, by \eqref{localextremumden} we obtain
	$\left( f_k(x_i+1)-f_k(x_i) \right) \cdot \left( f_k(x_i)-f_k(x_i-1) \right)<0 $, implying that $x_i$ is a local extremum.

	Now we move to  the situation  $\tan \theta(\gl_k, x_i)<0$. In this case, we have (cf. Figure \ref{fig:fixpoint}) 
\begin{equation*}
\ind_{I(\gl_k, c(x_i-1,x_i))}(\tan \theta(\gl_k, x_i))\;=\;0\,.
\end{equation*}		
	 As $\Xi^{(c)}(\lambda, b)>b$ for $b<0$, then
	$\Xi^{c(x_i-1,x_i)}(\gl_k, \tan \theta(\gl_k, x_i))>\tan \theta(\gl_k, x_i)$. Thus, using again \eqref{smaller_pi}, we infer
	\begin{equation*}
	i\pi -\frac{\pi}{2} \;<\;\theta (\gl_k, x_i) \;<\; i \pi\;<\; \theta(\gl_k, x_i+1)\;<\; i \pi +\frac{\pi}{2}\,.
	\end{equation*}
	Therefore, we obtain
	\begin{equation}\label{1branch}
	\begin{aligned}
	\tan \theta(\gl_k, x_i)\;&=\;- \frac{ c(x_i-1,x_i)\left( f_k(x_i)-f_k(x_i-1)\right)}{f_k(x_i-1)}\;<\;0\,,\\
	\tan \theta(\gl_k, x_i+1)\;&=\;- \frac{c(x_i,x_i+1)\left( f_k(x_i+1)-f_k(x_i)\right)}{f_k(x_i)}\;>\;0\,.\\
	\end{aligned}
	\end{equation}
	Assuming that $f_k(x_i-1) f_k(x_i)>0$, then we could multiply the two inequalities in \eqref{1branch} to obtain $\left( f_k(x_i)-f_k(x_i-1)\right) \cdot \left( f_k(x_i+1)-f_k(x_i) \right)<0$, leading to conclude that $x_k$ is a local extremum.

	Therefore, we only have to argue that $f_k(x_i-1) f_k(x_i)\leq 0$ cannot happen. The case $f_k(x_i-1) f_k(x_i)= 0$ cannot happen because $f_k(x_i-1)\neq 0$ and $\tan \theta(\gl_k, x_i)<0$.  Let us show that  $f_k(x_i-1) f_k(x_i)< 0$ cannot happen either. Suppose $f_k(x_i-1)<0$, and then $f_k(x_i)>0$. Hence
	\begin{equation*}
	\tan \theta(\gl_k,x_i)\;=\;-\frac{c(x_i-1,x_i) \left(f_k(x_i)-f_k(x_i-1)\right)}{f_k(x_i-1)}\;>\;0
	\end{equation*}
	which is a contradiction to $\tan \theta(\gl_k,x_i)<0$. For the other case  $f_k(x_i-1)>0$ and $f_k(x_i)<0$, the analysis is analogous.

	Now we deal with the case $\tan \theta(\gl_k, x_i)>0$. By the definition of $x_i$, we have that $(i-1)\pi< \theta(\gl_k,x_i)<
	i \pi-\frac{\pi}{2}$ and $
	i \pi< \theta(\gl_k,x_i+1)< i \pi+\frac{\pi}{2}$ since $\theta(\gl_k,x)$ and $\theta(\gl_k,x+1)$ are either in the same branch  or in adjacent branches of the tangent function by \eqref{def:varphi}.
	 Referring to Figure~\ref{fig:fixpoint} and the definition of $x_i$, we have
	\begin{equation}\label{3case:extremum}
	\begin{aligned}
	\tan \theta(\gl_k,x_i)\;&=\;- \frac{c(x_i-1,x_i) \left(f_k(x_i)-f_k(x_i-1) \right)}{f_k(x_i-1)}\;>\;c(x_i-1,x_i)\,,\\
	\tan \theta(\gl_k,x_i+1)\;&=\;- \frac{c(x_i,x_i+1) \left(f_k(x_i+1)-f_k(x_i) \right)}{f_k(x_i)}\;>\;0\,,
	\end{aligned}
	\end{equation}
	and thus
	\begin{equation*}
	\begin{gathered}
	f_k(x_i-1) f_k(x_i)\;<\;0\,,\\
	f_k(x_i) \left( f_k(x_i+1)-f_k(x_i) \right) \;<\;0\,.
	\end{gathered}
	\end{equation*}
	Hence, if $f_k(x_i)>0$, we have $f_k(x_i-1)<0$ and $f_k(x_i+1)-f_k(x_i)<0$ which implies $x_i$ is a local maximun.
	While if $f_k(x_i)<0$, we have $f_k(x_i-1)>0$ and $f_k(x_i+1)-f_k(x_i)>0$ which implies $x_i$ is a local minimum. We conclude the proof.
\end{proof}

\begin{lemma} \label{lema:nonlocalext}
	If $x \not \in \{ x_i\}_{i=1}^{k-1}$ (which includes the pair $\{ x_i-1,x_i\}$ when $\theta(\gl_k,x_i)=i \pi$) stated in Lemma~\ref{lema:extrema}, then $x$ is not a local extremum of the eigenfunction $f_k$.
\end{lemma}

\begin{proof}

	Let $x \in \lint 1, N \rint$ be as in the statement. Then there exists $i \in \bbN$ such that $(i-1) \pi < \theta(\gl_k, x)< \theta(\lambda_k, x+1) \le  i \pi$. If $\theta(\lambda_k, x+1) = i \pi$, by Lemma \ref{lema:extrema} both $x$ and $x+1$ would belong to $\{ x_i\}_{i=1}^k$, which contradicts the assumption. Thus we have
	$(i-1) \pi < \theta(\gl_k, x)< \theta(\lambda_k, x+1) < i \pi$, and thus $\tan \theta(\gl_k,x) \neq 0$, $\tan \theta(\gl_k,x+1) \neq 0$. 

	We first deal with the case $\tan  \theta(\gl_k,x)= \overline \infty$, i.e.\ $\theta(\gl_k,x)=(i-\frac{1}{2})\pi$. In this case, we have $f_k(x-1)=0$ and by the assumption of $x$,
	\begin{equation*}
	\tan \theta(\gl_k, x+1)\;=\;-\frac{c(x,x+1) \left(f_k(x+1)-f_k(x) \right)}{f_k(x)}\;<\;0\,.
	\end{equation*}
	Therefore, $\left( f_k(x+1)-f_k(x)\right) \cdot \left( f_k(x)-f_k(x-1)\right)>0$ and then $x$ is not a local extremum.

\medskip
	We move to the case where $\tan \theta(\gl_k, x)>0$ and $\tan \theta(\gl_k,x+1)$ satisfies one of the following conditions:  $\tan \theta(\gl_k,x+1)>0$ or $\tan \theta(\gl_k,x+1)=\overline \infty$.
	If $\tan \theta(\lambda_k, x)$ belongs to any of the intervals $(c,b_1)$, $[b_1, b_2]$ or $(b_2, \infty)$, referring to Figure \ref{fig:fixpoint}, the readers can check that it leads to a contradiction with the assumption   $\tan \theta(\gl_k,x+1)\in (0,\infty)\cup \{\overline{\infty}\}$ and the fact that $\theta(\lambda_k, x+1)-\theta(\lambda_k, x)<\pi/2$. Thus $\tan \theta(\lambda_k, x)\in (0,c]$, that is,
	\begin{equation*}
	\tan \theta(\gl_k,x)\;=\; -\frac{c(x-1,x) \left(f_k(x)-f_k(x-1) \right)}{f_k(x-1)}\;\leq \;c(x-1,x)
	\end{equation*}
	implying
	\begin{equation}\label{nonextre:cond1}
	f_k(x-1) f_k(x)\;\geq \;0\,.
	\end{equation}
	We do the analysis according to the value of $\tan \theta(\gl_k, x+1)$.\\
(1) If $\tan \theta(\gl_k, x+1)= \overline \infty$, then $f_k(x)=0$, $f_k(x-1) \neq 0$ and $f_k(x+1) \neq 0$  by \eqref{eq2extr}. Using $(\Delta^{(c)}f_k)(x)=-\gl_k f_k(x)$, we obtain
\begin{equation}
c(x-1,x)f_k(x-1)+c(x,x+1)f_k(x+1)=0
\end{equation}
which implies $x$ is not a local extremum. \\
(2)	If $\tan \theta(\gl_k, x)>0$ and 
	\begin{equation*}
	\tan \theta(\gl_k, x+1)\;=\;-\frac{c(x,x+1) \left( f_k(x+1)-f_k(x) \right)}{f_k(x)}\;>\;0\,,
	\end{equation*}
	we have
	\begin{equation}\label{nonextre:cond2}
	f_k(x) \left(f_k(x+1)-f_k(x) \right)\;<\;0\,,
	\end{equation}
	which in particular tells  $f_k(x)\neq 0$. Besides,
	since $\tan \theta(\lambda_k, x)\neq \overline{\infty}$,  we also infer $f_k(x-1)\neq 0$.	
	Now we show that $x$ is not a local extremum basing on \eqref{nonextre:cond1} and \eqref{nonextre:cond2}.
	If $f_k(x) > 0$, by \eqref{nonextre:cond1} we have  $f_k(x-1)> 0$.
	Then, due to $\tan \theta(\gl_k,x)>0$, we have  $f_k(x)-f_k(x-1)<0$. On the other hand, as $\tan \theta(\gl_k, x+1)>0$ and $f_k(x) > 0$, we have $f_k(x+1)-f_k(x)<0$. Therefore, $x$ is not a local extremum.
	If $f_k(x)< 0$, then by \eqref{nonextre:cond1} we have  $f_k(x-1) < 0$. Since $\tan \theta(\gl_k,x)>0$, it holds that $f_k(x)-f_k(x-1)>0$. As $\tan \theta(\gl_k, x+1)>0$ and $f_k(x)< 0$, we have $f_k(x+1)-f_k(x)>0$. Therefore, $x$ is not a local extremum.
	
	We now switch to the case $\tan \theta(\gl_k,x)>0$ and $\tan \theta(\gl_k,x+1)<0$. In this situation,  we have that (refer to Figure~\ref{fig:fixpoint})
	\begin{equation*}
	\tan \theta(\gl_k,x)\;=\; -\frac{c(x-1,x) \left(f_k(x)-f_k(x-1) \right)}{f_k(x-1)}\;>\;c(x-1,x)
	\end{equation*}
	and thus
	\begin{equation}\label{nonex:cond1}
	f_k(x-1)f_k(x)\;<\;0\,.
	\end{equation}
	On the other hand, since
	\begin{equation*}
	\tan \theta(\gl_k,x+1)\;=\;-\frac{c(x,x+1) \left(f_k(x+1)-f_k(x) \right)}{f_k(x)}\;<\;0\,,
	\end{equation*}
	we have
	\begin{equation}\label{nonex:cond2}
	f_k(x) \left(f_k(x+1)-f_k(x) \right)\;>\;0\,.
	\end{equation}
	If $f_k(x)>0$, by \eqref{nonex:cond1} we have $f_k(x-1)<0$. Additionally, by $\tan \theta(\gl_k,x)>0$ it holds that\break
	$f_k(x)-f_k(x-1)>0$. By \eqref{nonex:cond2} and $f_k(x)>0$, we have $f_k(x+1)-f_k(x)>0$. Therefore,  $x$ is not a local extremum.
	If $f_k(x)<0$, by \eqref{nonex:cond1} we have that  $f_k(x-1)>0$ and, by $\tan \theta(\gl_k,x)>0$, it holds that $f_k(x)-f_k(x-1)<0$. By $f_k(x)<0$  and $\tan \theta(\gl_k,x+1)<0$, we have $f_k(x+1)-f_k(x)<0$. Therefore,  $x$ is not a local extremum.\medskip

We are left to treat the case $\tan \theta(\gl_k,x)<0$.
		First of all, the case $\tan \theta(\gl_k,x)<0$ and $\tan \theta(\gl_k,x+1)=\overline \infty$ cannot occur, since  $\tan \theta(\gl_k,x+1)=\overline \infty$ implies $\tan \theta(\gl_k,x)=c(x-1,x)>0$, a contradiction.

		We move now to  the case
		\begin{equation*}
		\begin{aligned}
		\tan \theta(\gl_k,x)\;&=\;-\frac{c(x-1,x) \left(f_k(x)-f_k(x-1) \right)}{f_k(x-1)}\;<\;0\,,\\
		\tan \theta(\gl_k,x+1)\;&=\;-\frac{c(x,x+1) \left(f_k(x+1)-f_k(x) \right)}{f_k(x)}\;<\;0\,.\\
		\end{aligned}
		\end{equation*}
		If $f_k(x-1)f_k(x)>0$, by taking the product of the two inequalities above we have\break $\left(f_k(x)-f_k(x-1) \right) \cdot \left(f_k(x+1)-f_k(x) \right)>0$, and thus $x$ is not a local extremum.  We now argue that $f_k(x-1)f_k(x)<0$ cannot hold. Suppose that $f_k(x-1)f_k(x)<0$ holds.  If $f_k(x)<0$, then $f_k(x-1)>0$ and $\tan \theta(\gl_k,x)>0$, a contradiction to $\tan \theta(\gl_k,x)<0$. Else if $f_k(x)>0$, then we have $f_k(x-1)<0$ and $\tan \theta(\gl_k,x)>0$, a contradiction to $\tan \theta(\gl_k,x)<0$.

		We move to the case $\tan \theta(\gl_k, x)<0$ and $\tan \theta(\gl_k,x+1)>0$, which does not exist since $ (i-1) \pi < \theta(\gl_k,x)< \theta(\gl_k,x+1)< i\pi$.
		We conclude the proof.
\end{proof}

\begin{proof}[Proof of Proposition \ref{prop:eigfunc} concerning multi-monotonicity]

	Note that the number  of  monotonicity  equals to  the number of local extrema   plus one, according to our definition of local extremum stated above Proposition \ref{prop:eigfunc}.
	By  Lemma \ref{lema:numx}, Lemma \ref{lema:extrema}
	and Lemma \ref{lema:nonlocalext}, for the eigenvalue $\gl_j$, there are $j-1$ local extrema in the eigenfunction associated with $\gl_j$, and thus
	we conclude the proof.
\end{proof}

\section{The spectral gap under the assumption \eqref{LLN}} \label{sec:gapestimate}
In this section, our task is to prove Theorem \ref{th:gapshapeder}--\eqref{gap:asym}  under the assumption \eqref{LLN}.

\subsection{The first segment} We start dealing with those coordinates $x$ which are small.
\begin{proposition}\label{Prop:premierlem}
	If \eqref{LLN} holds, then for any $\gep \in (0,\pi/2)$ we have
	\begin{equation*}
	\limsup_{N\to \infty}\max_{u\in \left[0, \, (\frac{\pi}{2}-\gep)\ga^{-1/2}\right]} \; \left| B^{(N)}(\lceil u N     \rceil)- \sqrt{\alpha}\tan( \sqrt{\ga} u) \right|\;=\;0\,.
	\end{equation*}
\end{proposition}

\subsubsection{Denominator} We  first treat the denominator in the first term in the r.h.s.\ of \eqref{rel:scale}.
 Set
\begin{equation}\label{def:delta1}
\delta^{(1)}_N\;:=\;\max_{x\in \lint2, \,N\rint} \;\frac{r(x-1,x)}{N}\;\le\; \frac{1}{N}+ \frac{1}{N}\sup_{1\le n < m \le N} \left| (r(n,m)- (m-n)\right|
\end{equation}
which tends to zero by \eqref{LLN}.
If $z\in \left[ 0 , (\delta^{(1)}_N)^{1/2}\right]$,   setting $\delta^{(2)}_N=3 (\delta^{(1)}_N)^{1/2}$ we have
\begin{equation}\label{sandwichbound}
1+z\;\le\; \frac{1}{1-z} \;\le\; 1+ (1+\delta^{(2)}_N)z\,.
\end{equation}
Setting $(r')^{(N)}(x-1,x) \colonequals  (1+\delta^{(2)}_N)r^{(N)}(x-1,x)$ for $2 \le x \le N$,  in view of \eqref{sandwichbound} 
we define $\hat B^{(N)}$ and  $\tilde B^{(N)}$ by the following: (from now on we often drop the dependence in $N$ for ease of notations.)
\begin{equation}\label{def:Bhattilde}
\begin{cases}
\hat B(1)=\tilde B(1)=B(1)=0,\\
\hat B\left(x+1\right)\;=\; \hat B(x)+ N^{-1} r(x-1,x)\left(\hat B (x)\right)^2+\frac{\alpha}{N}\,, &\text{ for } x \ge 1\,, \\
\tilde B\left(x+1\right)\;=\; \tilde B(x)+ N^{-1} r'(x-1,x)\left(\tilde B (x)\right)^2+\frac{\alpha}{N}\,, &\text{ for } x \ge 1\,.
\end{cases}
\end{equation}
 Hence, as long as 
\begin{equation}\label{cond:sandwich}
 \max_{y\in \lint 0,x-1\rint} \tilde B(y)\le (\delta^{(1)}_N)^{-1/2}
\end{equation} 
  we have
\begin{equation}\label{sandbound:firstseg}
\hat B(x)\;\le\;   B(x)\;\le\; \tilde B(x)\,.
\end{equation}
Thus,  it is sufficient to prove Proposition
\ref{Prop:premierlem} concerning $\hat B$ and $\tilde B$.
Since $r'$ also satisfies the assumption~\eqref{LLN},  we just need to treat the case of  $\hat B$.

\subsubsection{Homogeneization}\label{subsec:homogen}
We provide a good approximation for $\hat B$, $\tilde B$  using the homogeneization method.

\begin{lemma}\label{lema:seclema}
	If \eqref{LLN} holds, then 
	for any $\gep>0$ we have
	\begin{equation}\label{1stseg:hom}
	\limsup_{N\to \infty} \, \max_{x\,:\, \frac{\sqrt{\alpha} r(1,x-1)}{N} \le   \frac{\pi}{2}-\gep }\,\left| \hat B(x)- \sqrt{\alpha}\tan\left( \frac{\sqrt{\alpha} r(1,x-1)}{N}\right) \right|\;=\;0\,.
	\end{equation}
\end{lemma}

\begin{proof}
	Setting 
	\begin{equation*}
	\begin{cases}
	Y(1)\; \colonequals\; 0\,, \\
	Y(x)\;\colonequals\;\sqrt{\alpha}\tan\left( \frac{\sqrt{\alpha} r(1,x-1)}{N}\right)\,, \; \forall \, x \ge 2 \,,
	\end{cases}
	\end{equation*}
	and using the formula for the tangent of the difference of two angles, we have that
	\begin{equation*}
	Y(x+1)\;=\;Y(x)+N^{-1}r(x-1,x)Y^2(x) + \frac{\alpha}{N}r(x-1,x)+  q_N(x)\,,
	\end{equation*}
	where $|q_N(x)|\le   C(\ga, \gep)  r(x-1,x)^2 N^{-2}$ for all $x$    stated in \eqref{1stseg:hom}.
	Set $w_N(1)=\gga(1)=0$, and for $x \ge 2$,
	\begin{equation*}
	\begin{aligned}
	w_N(x)\;&:=\; \frac{\alpha}{N}[x-1-r(1,x-1)]-\sum^{x-1}_{y=1}q_N(y)\,,\\
	\gamma(x)\;&:=\; \hat B(x)-Y(x)- w_N(x)\,.
	\end{aligned}
	\end{equation*}
 Then by \eqref{def:Bhattilde}	we have
	\begin{equation}\label{iterate:gamma-homo}
	\begin{aligned}
	\gamma(x+1)-\gamma(x)&\;=\;N^{-1}r(x-1,x) \left[\hat B(x)^2 - Y(x)^2\right]\\
	&\;=\;N^{-1}r(x-1,x) \cdot \left[2Y(x)+\gamma(x)+w_N(x)\right]\cdot \left[\gamma(x)+w_N(x)\right]\,.
	\end{aligned}
	\end{equation}
	 In our range of $x$, we know that the sequence $Y(x)$ is uniformly bounded and  $|w_N(x)|\le \delta_N
		\colonequals  \ga \delta_N^{(0)}+ 2 C(\ga,\gep)\delta_N^{(1)}$,
	where $\delta_N$ tends to zero by \eqref{LLN} and \eqref{def:delta1}.

	Now we argue that
	$|\gamma(x)|\le 1$ for all $x$ stated in \eqref{1stseg:hom} by induction. By definition, it holds for $x=1,2$ (for all $N$ big enough). Now we suppose it holds for all $k \le x$.  Then  
	by \eqref{iterate:gamma-homo}	 and  $ \left| 2Y(x)+\gamma(x)+w_N(x)\right| \le   C=C(\ga,\gep)$, 
	we have for all $k \le x$
		\begin{equation*}
		|\gamma(k+1)|\;\le\; \left(1+ C N^{-1}r(k-1,k)\right) |\gamma(k)|+C N^{-1}r(k-1,k)\delta_N\,.
		\end{equation*}
		Iterating this inequality from $k=x$ backward to $k=2$, we obtain
		\begin{align*}
		\vert \gga(x+1) \vert
		\;&\le\;
		\vert \gga(2) \vert \prod_{j=2}^x \left(1+\frac{C r(j-1,j)}{N} \right)+ \sum_{j=2}^x \frac{C r(j-1,j)}{N} \delta_N \prod_{i=j+1}^x \left(1+ \frac{C r(i-1,i) }{N}\right)\\
		\;&\le\;
		\vert \gga(2) \vert \exp \left( \ \frac{C r(1,x)}{N}  \right)
		+  \frac{Cr(1,x)}{N} \delta_N \exp \left( \frac{C r(1,x)}{N} \right) \;\le\; C' \delta_N\,.
		\end{align*}
		Therefore, the assumption $|\gamma(x+1)|\le 1$ is verified, and we can conclude that
		\begin{equation*}
		\limsup_{N\to \infty} \max_{x: \frac{\sqrt{\alpha} r(1,x-1)}{N} \le \pi/2-\gep } |\hat B(x)-Y(x)|
\;\le \;
\limsup_{N\to \infty} \max_{x: \frac{\sqrt{\alpha} r(1,x-1)}{N} \le \pi/2-\gep }	\vert \gga(x)\vert	+ \vert w_N(x) \vert
		\;=\;0\,.
		\end{equation*}
\end{proof}
 
 \begin{proof}[Proof of Proposition \ref{Prop:premierlem}]
 Since $r'$ satisfies \eqref{LLN}, we have 
   a version of Lemma \ref{lema:seclema} concerning $\tilde B$ with  $r$ replaced by $r'$, i.e.,
  \begin{equation}\label{1stseg:homrprime}
	\limsup_{N\to \infty}\,\max_{x\,:\, \frac{\sqrt{\alpha} r'(1,x-1)}{N} \le   \frac{\pi}{2}-\gep } \, \left| \tilde B(x)- \sqrt{\alpha}\tan\left( \frac{\sqrt{\alpha} r'(1,x-1)}{N}\right) \right|\;=\;0\,,
	\end{equation}
 and then  the condition \eqref{cond:sandwich} is verified. 
 Since
\begin{equation*}
  \left\vert \tan \left(t\sqrt{\ga}\right)-\tan\left(s\sqrt{\ga}\right) \right\vert \;\le\; C(\gep) \sqrt{\ga} \left\vert t-s \right\vert  \quad \text{ for } 0 \;\le\;  s\sqrt{\ga} \,, t\sqrt{\ga} \;\le\; \frac{\pi}{2}-\gep\,, 
\end{equation*}  
  we conclude the proof by Lemma \ref{lema:seclema}, \eqref{1stseg:homrprime}, \eqref{LLN}, $r'=(1+\delta_N^{(2)}) r$  and \eqref{sandbound:firstseg}.

 \end{proof}

 \subsection{The second segment}
To deal with those $x$ which are not covered in the interval stated in \eqref{1stseg:hom}, we set 
\begin{equation}\label{def:AinverB}
A^{(N)}(x)\;\colonequals\; \frac{1}{B^{(N)}(x)}  \,.
\end{equation}
 The motivation of this mapping 
is to  avoid dealing with the explosion of the tangent function at the neighbor of $\pi/2$. We often  drop the superscript to lighten the notations. By \eqref{rel:scale} we have

\begin{equation}\label{iteration:A}
N \left[ A\left(x+1 \right)-A\left( x\right) \right]\;=\;
\frac{\ga r(x-1,x) A\left(x \right)N^{-1}-r(x-1,x)-\ga A\left( x\right)^2}
{1+\ga A\left(x\right)N^{-1}-\ga r(x-1,x)N^{-2}}\,.
\end{equation}
 Define now
\begin{equation}\label{def:taus}
\begin{aligned}
\tau_1 \;&\colonequals\; \inf \left\{x \in \lint 1, N \rint: \;    \frac{\sqrt{\alpha} r(1,x)}{N} \ge \frac{ \pi}{4} \right\}\,,\\
\tau_2 \;&\colonequals\; \sup \left\{x \in \lint 1, N \rint: \;    \frac{\sqrt{\alpha} r(1,x)}{N} \le \frac{3 \pi}{4} \right\}\,,\\
\tau \;&\colonequals\; \inf \left\{x \in \lint 1, N \rint: \;    A(x) \le -2 \ga^{-1/2} \right\}-1\,.\\
\end{aligned}
\end{equation}
For $x \in \lint \tau_1, \tau\rint$, by \eqref{def:delta1} we have
\begin{equation}\label{mondecrease:A}
1+ \frac{\ga A\left( x \right)}{N}-\frac{\ga r(x-1,x)}{N^2}
\;\ge\;
1- \frac{2 \sqrt \ga}{N}-\frac{\ga r(x-1,x)}{N^2}\;>\;0\,.
\end{equation}
For all $x \in \lint \tau_1,\tau\rint$, by induction we know that $\ga A(x)N^{-1} <1$ for all $N$ sufficiently large, and thus
$A( x)_{\tau_1 \le x \le  \tau}$ is monotone decreasing. By Proposition \ref{Prop:premierlem},  
we have for all $\tau_1 \le x \le  \tau$
\begin{equation*}
\vert A(x)\vert \;\le\; \max(2\ga^{-1/2},2) \;\equalscolon\; K\,.
\end{equation*}
With the same spirit as in \eqref{def:Bhattilde}, 
in view of \eqref{iteration:A} we define
\begin{equation}\label{sandw:As}
\begin{cases}
\hat A(\tau_1)\;=\;\tilde A(\tau_1)\;=\;A(\tau_1)\,,\\
\hat A\left(x+1\right)\; \colonequals  \;  \hat A(x)- N^{-1} \hat \ga \left(\hat A (x)\right)^2-\frac{\hat r(x-1,x)}{N}\,, \quad x \in \lint \tau_1, \, \tau \wedge \tau_2 \rint\,,\\
\tilde A\left(x+1\right)\;  \colonequals  \;  \tilde A(x)- N^{-1} \tilde \ga \left(\tilde A (x)\right)^2-\frac{\tilde r(x-1,x)}{N}\,, \quad x \in \lint \tau_1, \, \tau \wedge \tau_2 \rint\,.
\end{cases}
\end{equation}
where (we drop the dependence on $N$ to ease the notations)
\begin{equation}\label{def:hattildegar}
\begin{aligned}
\hat \ga \;&\colonequals \frac{\ga}{1-K\ga N^{-1}-\ga \delta_N^{(1)}N^{-1}}\,, \quad \hat r(x-1,x)\;\colonequals\; \frac{ r(x-1,x)}{1-K\ga N^{-1}-\ga \delta_N^{(1)}N^{-1}} \left( 1+\frac{\ga K}{N}\right)\,,\\
\tilde \ga \;&\colonequals\; \frac{\ga}{1+K\ga N^{-1}}\,, \quad \quad \quad \quad \quad  \quad \tilde r(x-1,x)\;\colonequals\; \frac{ r(x-1,x)}{1+K\ga N^{-1}} \left( 1-\frac{\ga K}{N}\right)\,.
\end{aligned}
\end{equation}
Moreover, by \eqref{sandw:As} we have
\begin{align}\label{compare:AhatA}
A(x+1)-\hat A(x+1)
\;&\ge\; \left[A(x)- \frac{\hat \ga}{N} A(x)^2-\frac{\hat r(x-1,x)}{N}\right]-\left[ \hat A(x)-\frac{ \hat a}{N} \hat A(x)^2- \frac{\hat r(x-1,x)}{N}  \right] \nonumber\\ 
\;&=\;\left[ A(x)-\hat A(x) \right]\cdot \left[1- \frac{\hat \ga}{N} \left( A(x)+\hat A(x)\right) \right]\,,
\end{align}
and then we can argue by induction that the r.h.s.\ above is nonnegative by the fact that  $A\left( x\right)_{\tau_1 \le x \le  \tau}$ and $\hat A\left( x\right)_{\tau_1 \le x \le  \tau}$ are  decreasing.
By a similar version  of \eqref{compare:AhatA}  between $A$ and $\tilde A$, 
we have
\begin{equation}\label{sand:hatAtilde}
\hat A(x) \;\le\; A(x) \;\le\; \tilde A(x)\, \quad \forall \, x \in \lint \tau_1,\, \tau \wedge \tau_2 \rint\,.
\end{equation}
Note that $\hat r$ and $\tilde r$ also satisfy the assumption \eqref{LLN}, and $\lim_{N \to \infty} \hat \ga=\lim_{N \to \infty} \tilde \ga=\ga$. It is sufficient to treat $\hat A$.

\begin{lemma}\label{lema:2seghom}
	If \eqref{LLN} holds, for $\sqrt{\ga} \ge \frac{3\pi}{4}$ we have
	\begin{equation}\label{2nd:hom}
	\limsup_{N \to \infty}\; \max_{x \in \lint \tau_1, \, \tau_2 \wedge \tau \rint}  \; \left|\hat A(x)- \ga^{-1/2} \tan \left(\frac{\pi}{2}- \frac{(x-1)\sqrt{ \ga} }{N} \right) \right|\;=\;0\,.
	\end{equation}
\end{lemma}

\begin{proof}
	Note that under the  assumptions \eqref{LLN} and $\sqrt{\ga} \ge \frac{3\pi}{4}$, we have
	\begin{equation*}
	\limsup_{N \to \infty} \; \max_{x \in \lint \tau_1, \, \tau_2 \wedge \tau\rint } \; \left \vert
	\ga^{-1/2} \tan \left(\frac{\pi}{2}- \frac{(x-1)\sqrt{ \ga} }{N} \right)
	-
	\hat\ga^{-1/2} \tan \left(\frac{\pi}{2}- \frac{(x-1)\sqrt{\hat \ga} }{N} \right)
	\right\vert\;=\;0\,.
	\end{equation*}
	Thus it is sufficient to prove \eqref{2nd:hom} with $\ga, \; r$ replaced by $\hat \ga, \; \hat r$ respectively.

	Setting
	\begin{equation*}
	Z(x)\; \colonequals \;\hat \alpha^{-1/2}\tan\left(\frac{\pi}{2}- \frac{(x-1)\sqrt{\hat \alpha}  }{N}\right)  \quad \forall \, x \in \lint \tau_1, \, \tau_2 \wedge \tau \rint\,,
	\end{equation*}
	we have
	\begin{equation}\label{recur:Z}
	Z(x+1)\;=\; Z(x)-N^{-1}- \frac{\hat \alpha}{N} Z(x)^2+  \hat q_N(x)
	\end{equation}
	where $|\hat q_N(x)|\le { \bar C(\ga)}  N^{-2}$. 
For $x \in \lint \tau_1, \, \tau_2 \wedge \tau \rint$,	setting
	\begin{equation*}
	\begin{split}
	\hat w_N(x)& \;:=\; \sum_{y=\tau_1}^{x-1}\frac{1}{N}[\hat r(y-1,y)-1]+\sum^{x-1}_{y=\tau_1} \hat q_N(y)\,, \\
	\Gamma(x)\;&:=\; \hat A(x)-Z(x)+ \hat w_N(x)\,,
	\end{split}
	\end{equation*}
where $\hat w_N(\tau_1) \colonequals 0$,   
by \eqref{sandw:As} and \eqref{recur:Z}	we have
	\begin{equation}\label{iterate:Gamma}
	\begin{split}
	\Gamma(x+1)-\Gamma(x)& \;=\; -\hat \ga N^{-1} \left[ \hat A(x)^2-Z(x)^2 \right]
	\\
	&\;=\;-\hat \ga N^{-1}  \left[2Z(x)+\Gamma(x)-\hat w_N(x)\right ]\cdot \left[\Gamma(x)- \hat w_N(x)\right]\,.
	\end{split}
	\end{equation}
	 Moreover, by the definition of $\hat r$ in \eqref{def:hattildegar} and \eqref{LLN}, we have
	\begin{equation}\label{error:Gamma}
\begin{aligned}
\vert \hat w_N(x) \vert &\; \le\; \bar C(\ga)N^{-1}+
\max_{x \in \lint \tau_1+1, \, N+1\rint }\left\vert \frac{1}{N} \sum_{y=\tau_1}^{x-1} \left[\hat r(y-1,y)-1\right] \right\vert \\
&\;\le\; 
\bar C(\ga)N^{-1}+ \max_{x \in \lint \tau_1+1, \, N+1\rint } \frac{1}{N} \left[ \vert \hat r(1,x-1)-(x-2) \vert + \vert \hat r(1,\tau_1-1)-(\tau_1-2)\vert\right]  \\
&\;\le\;
\bar C(\ga)N^{-1}+ 2 \delta_N^{(0)}  \left(1+\frac{\ga K}{N} \right) \left( 1-K\ga N^{-1}-\ga \delta_N^{(1)}N^{-1}\right)^{-1} \\
 & \quad + 2 \left[\left(1+\frac{\ga K}{N} \right) \left( 1-K\ga N^{-1}-\ga \delta_N^{(1)}N^{-1}\right)^{-1}-1\right] \equalscolon \hat \delta_N\,.
\end{aligned}	
\end{equation}
Recalling $\delta_N^{(0)}$ and $\delta_N^{(1)}$ defined in \eqref{LLN}
and \eqref{def:delta1} respectively, we know that $\hat \delta_N$ tends to $0$.
	 Now we argue that 
	$|\Gamma(x)|\le 1$ for all $x$ stated in \eqref{2nd:hom} by induction. 
Note that
		\begin{equation}\label{sm:Gamtau1}
		\Gamma(\tau_1)\;=\;\hat A(\tau_1)-Z(\tau_1)\;=\; \hat B(\tau_1)^{-1}-\hat \ga^{-1/2}\left[ \tan \left( \frac{(\tau_1-1) \sqrt{ \hat \ga}}{N}\right)\right]^{-1}
		\end{equation}
		which tends to $0$ by Lemma \ref{lema:seclema} and \eqref{LLN}, implying $\vert \Gamma(\tau_1) \vert \le 1$.
			 Now we suppose $|\Gamma(x)|\le 1$ for all $k \in \lint \tau_1, x\rint$.  Note that
under the assumption \eqref{LLN},  in our range of $x$,
	$Z(x)$ is uniformly bounded.	
	Then by \eqref{iterate:Gamma} and   \eqref{error:Gamma},  for some constant $C=C(\ga)$ 
		we have for all $k \in \lint \tau_1, x \rint$
		\begin{equation*}
		|\Gamma(k+1)|\;\le\; \left(1+ C \hat \ga  N^{-1} \right) |\Gamma(k)|+C \hat \ga N^{-1}  \hat \delta_N\,.
		\end{equation*}		
		Iterating this inequality from $k=x$ backward to $k=\tau_1$, we obtain
\begin{equation}	
\begin{aligned}
		\vert \Gamma(x+1) \vert
		& \; \le\;
		\vert \Gamma(\tau_1) \vert \prod_{j=\tau_1}^x \left(1+\frac{C \hat \ga}{N} \right)
		+ \sum_{j=\tau_1}^x \frac{C \hat \ga}{N} \hat \delta_N \prod_{i=j+1}^x \left(1+ \frac{C \hat \ga }{N}\right)\\
		&\;\le\;
		\vert\Gamma(\tau_1)\vert \exp \left( C \hat \ga  \right)
		+  C \hat \ga \hat \delta_N \exp \left( C \hat \ga \right)\,.
\end{aligned}
\end{equation}	
By \eqref{sm:Gamtau1}, 		
		the assumption on $|\Gamma(x+1)|\le 1$ is verified, and thus  we obtain
\begin{equation*}
\limsup_{N\to \infty}\; \max_{x \in \lint \tau_1,\, \tau_2 \wedge \tau \rint}\, \left|\Gamma(x)\right|\;=\;0\,,	
\end{equation*}		
	which allows us to conclude the proof by	 \eqref{error:Gamma} and $\hat A(x)-Z(x)=\Gamma(x)-\hat w_N(x)$.
\end{proof}

Note that a version of Lemma \ref{lema:2seghom} concerning $\tilde A$, defined in \eqref{sandw:As}, holds by the same argument. Furthermore,   by \eqref{sand:hatAtilde},  a version of Lemma \ref{lema:2seghom} concerning $ A$ also holds, and thus
 $\tau_2< \tau$ where we assume $ \sqrt{\ga}> \frac{3\pi}{4}$. That is to say, we have:
\begin{lemma}\label{lema:2seghomA}
	Assuming \eqref{LLN}, for $  \sqrt{\ga} > \frac{3\pi}{4}$  the following holds:
	\begin{equation}
	\limsup_{N \to \infty}\; \max_{x \in \lint  \tau_1, \,\tau_2\rint}  \; \left| A(x)- \ga^{-1/2} \tan \left(\frac{\pi}{2}- \frac{(x-1)\sqrt{ \ga} }{N} \right)  \right|\;=\;0\,.
	\end{equation}
\end{lemma}

By Lemma \ref{lema:2seghomA}, if $\sqrt{\ga}> \frac{3 \pi}{4}$, then we have $\tau_2<N$ and
\begin{equation*}
\lim_{N \to \infty} B(\tau_2)\;=\;-\sqrt{\ga}\,.
\end{equation*}

 \subsection{The last segment}
From now on, we move to treat $x \in \lint \tau_2, N+1\rint$.
We set
\begin{equation*}
\tau' \colonequals \sup \left\{x \in \lint \tau_2,\, N+1\rint : \; B(x)\le 2 \sqrt{\ga}\right\}\,.
\end{equation*}
Combining \eqref{def:delta1} with the fact that $B(x) \le 2 \sqrt{\ga}$ for $x \in \lint \tau_2, \tau'\rint$, by \eqref{resscale:B} we know that $B$ is  increasing in $\lint \tau_2, \tau'\rint$, and thus
\begin{equation}
\vert B(x)\vert \;\le\; K' \colonequals 2 \sqrt{\ga}\,, \quad \forall \, x \in \lint \tau_2, \, \tau' \rint\,.
\end{equation}
With the sandwich trick as in \eqref{sandw:As}, an analogous proof as that in Lemma \ref{lema:seclema}, and Lemma \ref{lema:2seghomA} for $x=\tau_2$, we have the following lemma.
\begin{lemma}\label{lema:3seghomB}
	If \eqref{LLN} holds and  $  \sqrt{\ga}> \frac{3\pi}{4}$,  we have
	\begin{equation}
	\limsup_{N \to \infty}\; \max_{x \in \lint  \tau_2, \, \tau' \rint}  \; \left| B(x)- \ga^{1/2} \tan \left(\frac{\sqrt{ \ga } r(1,x-1)}{N} \right) \right|\;=\;0\,.
	\end{equation}
\end{lemma}
 
 With all the ingredients above, we are ready to prove Theorem \ref{th:gapshapeder}--\eqref{gap:asym}. 
\begin{proof}[Proof of   Theorem \ref{th:gapshapeder}--\eqref{gap:asym}]
	With the assumption stated in Lemma \ref{lema:3seghomB} 
	we know that $N+1= \tau'$.
	Set
	\begin{equation*}
	x_0 \colonequals \inf \left\{  x \in \lint 1, N\rint: \ B(x)>0, \, B(x+1) \le 0 \right\}\,,
	\end{equation*}
	which exists by Lemma \ref{lema:2seghomA} and $B(x)=1/A(x)$.
	Note that $B$ is strictly increasing in $\lint 1,\, x_0 \rint$ by  \eqref{resscale:B} and the two lines below \eqref{mondecrease:A} since $B(x)=1/A(x)$.  Moreover, $B$ is 
	 strictly increasing in $\lint x_0+1,\, N+1\rint$ due to \eqref{resscale:B} and  $B(N+1) \le 2 \sqrt{\ga}$ in this interval by Lemma \ref{lema:3seghomB}.

	To obtain an estimate on the spectral gap $\gl_1=\ga/N^2$, we recall that $B(x)=b(\gl,x)N$, and that $b(\gl, x)$ and $\theta(\gl,x)$ are related by \eqref{def:Xi} and \eqref{def:maptheta}.
	By Lemma \ref{lema:thetacont} and Lemma \ref{lema:3seghomB}, for any given $\gep>0$, for all $N$ sufficiently large we have
	\begin{equation*}
	\left\vert \sqrt{\ga}-\pi \right\vert \; \le \;\gep\,,
	\end{equation*}
	which allows  us to obtain \eqref{gap:asym}.

\end{proof}

\section{The  principal eigenfunction and its  derivative  under the assumption~\eqref{LLN}}\label{sec:eigshapeder}
 In this section, our goal is to prove  \eqref{eigfun:shape} and
\eqref{approx:der2eigfuns} in Theorem \ref{th:gapshapeder}. To lighten notations, 
let $g=g_1^{(N)}$ be the principal eigenfunction corresponding to the spectral gap as stated in Theorem \ref{th:gapshapeder}.
Note that the function
\begin{equation}\label{fun:h}
h(x)\;=\;h_N(x) \;=\; \cos \left( \frac{\pi(x-1/2)}{N} \right)
\end{equation}
is the principal eigenfunction corresponding to the spectral gap when $r(j-1,j) \equiv 1$, and in this case  the spectral gap is
\begin{equation*}
 \overline \gl  \;\colonequals\; 2\left[1-\cos\left( \frac{\pi}{N}\right) \right]
\;=\;\frac{\pi^2}{ N^2}+O \left( \frac{1}{N^4} \right)\,.
\end{equation*}
Moreover, by \eqref{def:blambda} and $b(\gl,x)=B(x)/N$, for $x \ge 2$ we have
\begin{equation}\label{rec:randeigfun1st}
g(x) \;=\; \left[1- r(x-1,x)N^{-1} B \left( x
\right) \right] g(x-1)\,,
\end{equation}
where we recall $g(1)=1$.
We define $(\overline B(x))_{1 \le x \le N+1}$ and $(\overline b(\overline \gl, x))_{1\le x \le N+1}$ with $\overline b(\overline \gl, x)=\overline B(x)/N$ 
by setting $\overline B \left( 1 \right)=0$ and for $x \in \lint 1, N\rint$ (with the convention $1/0 \colonequals \overline \infty$), 
\begin{equation}\label{def:breveB}
\begin{aligned}
\overline b(\overline \gl, x ) \;&\colonequals \;-\frac{h(x)-h(x-1)}{h(x-1)} \,,\\
\overline B\left(x+1\right) \;&\colonequals\;
\frac{  \overline B\left(x\right)}{1-\overline B \left(x\right)N^{-1 }}+ N \overline \gl\,,
\end{aligned}
\end{equation}
which are 
the corresponding deterministic versions of \eqref{def:blambda} and \eqref{rel:scale}  with $\gl$ replaced  by $\overline \gl$.
Similar to \eqref{rec:randeigfun1st}, we have
\begin{equation}\label{rec:detereigfun1st}
h(x)\;=\; \left[1- N^{-1} \overline B \left( x \right) \right] h(x-1)\,.
\end{equation}

\subsection{ The shape and derivative in the first segment}
\begin{lemma}\label{lema:1stapproxfun}
 If \eqref{LLN} holds,  we have
\begin{equation}\label{autofuns:approx}
\begin{gathered}
\lim_{N \to \infty} \; \max_{x \in \lint 1,\, \tau_1 \rint} \; \left\vert g(x)- h(x) \right\vert \;=\;0\,,\\
\end{gathered}
\end{equation}
and
\begin{equation}\label{der:approx}
\lim_{N \to \infty} \; \max_{x \in \lint 1, \, \tau_1 \rint} \; \left\vert N (c\nabla g)(x)- N(\nabla h)(x) \right\vert\;=\;0\,.
\end{equation}
where $\tau_1$ is defined in \eqref{def:taus} with $\ga/N^2$ being the spectral gap of the operator $\Delta^{(c)}$.
\end{lemma}

\begin{proof}
We first deal with \eqref{autofuns:approx}.
Combining \eqref{rec:randeigfun1st} with \eqref{rec:detereigfun1st} and
writing $u(x)\colonequals h(x)-g(x)$, we have
\begin{equation}\label{iterate:ux}
u(x)
\;=\;
u(x-1) \left[1-\frac{r(x-1,x)B \left( x\right)}{N}\right]+\frac{h(x-1)}{N} \left[ r(x-1,x)B \left(x \right)-  \overline B \left(x \right)\right]\,.
\end{equation}
Iterating the equation above, we obtain
\begin{equation}
\begin{aligned}\label{backward:difgapfun}
u(x)&=u(1)\prod_{k=1}^{x-1} \left[ 1-\frac{r(k,k+1)B\left(k+1 \right)}{N}\right]\\
& \quad+
\sum_{k=1}^{x-1} \frac{h(k)}{N} \left[ r(k,k+1)B \left(k+1 \right)- \overline B \left( k+1
\right) \right] \prod_{j=k+1}^{x-1} \left[ 1- \frac{r(j,j+1)B \left( j+1\right)}{N} \right].
\end{aligned}
\end{equation}
Concerning the first term in the r.h.s.\ of \eqref{backward:difgapfun}, 
by \eqref{LLN}, Theorem \ref{th:gapshapeder} and Proposition \ref{Prop:premierlem}
for all $N$ sufficiently large
we have  
\begin{equation}\label{bound:B1seg}
0 \le  B\left( k\right) \le 2 \sqrt \ga\,, \quad \forall \, k \in \lint 1, \tau_1\rint\,,
\end{equation}
and then
\begin{equation*}
\left \vert u(1)\prod_{k=1}^{x-1} \left[ 1-\frac{r(k,k+1)B\left(\frac{k+1}{N} \right)}{N}\right] \right \vert \;\le\; \vert u(1)\vert\;=\; 1-\cos \left( \frac{\pi}{2N}\right) \;\le\; \frac{C}{N^2}\,.
\end{equation*}
The second term in the r.h.s.\ of \eqref{backward:difgapfun} is
\begin{equation}\label{def:sumremaining}
\sum_{k=1}^{x-1}\frac{h(k)}{N}
\left\{
\left[
r(k,k+1)-1
\right] B \left( k+1 \right)
+
\left[ B \left( k+1 \right)- \overline B \left( k+1 \right)\right]
\right\}v_{k+1}\,,
\end{equation}
with
\begin{equation*}
v_k \;\colonequals \; \prod_{j=k}^{x-1} \left[ 1- r(j,j+1)B \left( j+1\right)N^{-1} \right]
\end{equation*}
where we ignore the dependence in $x$ and $ N$ to lighten notations. We first deal with the term $\sum_{k=1}^{x-1}\frac{h(k)}{N}
\left[ B \left( k+1 \right)- \overline B \left( k+1 \right)\right]
v_{k+1}$ in \eqref{def:sumremaining}.
Note that by Proposition \ref{Prop:premierlem} and \eqref{def:delta1},  we have
\begin{equation}\label{bound:vk}
 \vert v_k \vert \;\le\; 1\,, \quad \forall\, k \in \lint 1,\, \tau_1\rint\,.
\end{equation}
Since $\vert h(k) \vert \le 1$, $\vert v_k \vert \le 1$ and by Lemma \ref{lema:seclema}
\begin{equation*}
\limsup_{N \to \infty}\; \max_{x \in \lint 1, \,\tau_1\rint}\; \left \vert B(x)-\overline B(x) \right\vert\;=\;0\,,
\end{equation*}
then we have
\begin{equation*}
\limsup_{N \to \infty}\; \max_{x \in \lint 1,\, \tau_1\rint} \; \left\vert \sum_{k=1}^{x-1}\frac{h(k)}{N}
\left[ B \left( k+1 \right)- \overline B \left( k+1 \right) \right]
v_{k+1} \right\vert\;=\;0\,.
\end{equation*}
Now we turn to the other term in \eqref{def:sumremaining}. Setting $w_n \colonequals \sum_{k=1}^n [r(k,k+1)-1]$ with $w_0 \colonequals 0$, we have
\begin{equation}\label{eigfun:shapemart}
\begin{aligned}
&\sum_{k=1}^{x-1}\frac{h(k)}{N}
\left[
r(k,k+1)-1
\right] B \left( k+1 \right)
v_{k+1}
\;=\;
\sum_{k=1}^{x-1}\frac{h(k)}{N}
B \left( k+1 \right)
v_{k+1} (w_k-w_{k-1})\\
=\;&
\frac{1}{N} \sum_{k=1}^{x-1}  \left[ h(k)
B \left( k+1 \right)
v_{k+1}- h(k+1)
B \left( k+2 \right)
v_{k+2} \right]w_k+ \frac{h(x)}{N}  B \left( x+1\right) v_{x+1}w_{x-1}\,.
\end{aligned}
\end{equation}
To deal with the first term in the r.h.s.\ of \eqref{eigfun:shapemart}, 
note that by triangle inequality, 
\begin{equation}\label{3terms:inserted}
\begin{aligned}
&\left \vert  h(k)
B \left( k+1 \right)
v_{k+1}-
h(k+1)
B \left( k+2 \right)
v_{k+2}
\right\vert
\; \le\; 
\left \vert  h(k)
B \left( k+1 \right)
v_{k+1}-
h(k)
B \left( k+1 \right)
v_{k+2}
\right\vert\\
&\, +
\left \vert
h(k)
B \left( k+1 \right)
v_{k+2}-
h(k)
B \left( k+2 \right)
v_{k+2}
\right\vert +\left \vert
h(k)
B \left( k+2 \right)
v_{k+2}-
h(k+1)
B \left( k+2 \right)
v_{k+2}
\right\vert \\
 & \, \le
2 \sqrt{\ga} \left \vert  v_{k+1}-v_{k+2} \right\vert +
\left\vert  B \left(k+1 \right) -B \left(k+2 \right)\right\vert+ 2 \sqrt{\ga} \left \vert h(k)-h(k+1) \right \vert\,,
\end{aligned}
\end{equation}
where we have used \eqref{bound:vk} and \eqref{bound:B1seg}.
 Moreover, we have
\begin{equation}\label{3erros}
\begin{aligned}
\vert v_{k+1}-v_{k+2} \vert \; &\le\; r(k+1,k+2)N^{-1}B(k+2)\,,\\
\left\vert  B \left(k+1 \right) -B \left(k+2 \right)\right\vert \; &\le \;  \frac{\ga}{N}+\frac{N^{-1}r(k,k+1) B(k+1)^2}{1-N^{-1}r(k,k+1)B(x)}\,,\\
\left \vert h(k)-h(k+1) \right \vert \;&\le\; \pi N^{-1}\,.
\end{aligned}
\end{equation}
By \eqref{3terms:inserted} and \eqref{3erros}, the summation term in the r.h.s.\ of \eqref{eigfun:shapemart}  can be bounded from above by
\begin{equation}\label{bounds:sumterms}
\begin{aligned}
&\frac{1}{N} \sum_{k=1}^{x-1} \left[ 2 \sqrt{\ga} \left \vert  v_{k+1}-v_{k+2} \right\vert +
\left\vert  B \left(k+1 \right) -B \left(k+2 \right)\right\vert+ 2 \sqrt{\ga} \left \vert h(k)-h(k+1) \right \vert \right] \cdot\vert w_k \vert \\
\le \; &
\frac{2 \sqrt{\ga} }{N} \sum_{k=1}^{x-1} r(k+1,k+2)N^{-1} B(k+2) \vert  w_k\vert
 + \frac{\ga}{N^2} \sum_{k=1}^{x-1} \vert w_k \vert \\
 & + \frac{1}{N}\sum_{k=1}^{x-1} \vert w_k \vert  \frac{N^{-1}r(k,k+1) B(k+1)^2}{1-N^{-1}r(k,k+1)B(x)}
 + \frac{2 \sqrt \ga \pi}{N^2} \sum_{k=1}^{x-1} \vert w_k \vert \,.
\end{aligned}
\end{equation} 
 Concerning the first summation term above,  by \eqref{bound:B1seg} and \eqref{LLN} we have 
\begin{equation*}
\frac{1}{N} \sum_{k=1}^{x-1} r(k+1,k+2)N^{-1} B(k+2) \vert  w_k\vert 
\;\le\; \delta_N^{(0)} 2 \sqrt{\ga} \frac{1}{N} \sum_{k=1}^{x-1}r(k+1,k+2) \;\le\; 4 \sqrt{\ga} \delta_N^{(0)}
\end{equation*}
which tends to zero. Both the second and fourth summation terms in \eqref{bounds:sumterms} tend to zero, which is a consequence of \eqref{LLN}. Concerning the third summation term in \eqref{bounds:sumterms}, by \eqref{bound:B1seg} and \eqref{def:delta1} it can be bounded from above by
\begin{equation*}
\max_{1 \le j \le x-1} \frac{\vert w_j\vert}{N}  \sum_{k=1}^{x-1}  C(\ga) \frac{r(k,k+1)}{N} \;\le \; 2 C(\ga)  \max_{1 \le j \le x-1} \frac{\vert w_j\vert}{N} \; \le  \; 2 C(\ga) \delta_N^{(0)}\,.
\end{equation*}
Additionally, the last term in the r.h.s.\ of \eqref{eigfun:shapemart} also tends to zero due to $\vert h(x+1) \vert \le 1$, \eqref{bound:B1seg}, \eqref{bound:vk} and \eqref{LLN}.
 Therefore,  we conclude the proof of \eqref{autofuns:approx}.

\medskip

Now we move to deal with \eqref{der:approx}.
By \eqref{grad:g} we have
\begin{equation}
\label{diff:derivatives}
\begin{aligned}
\left \vert N (c \nabla g)(x)-N (\nabla h)(x)  \right \vert 
\;=\;&
\left \vert -N\gl_1 \sum_{k=1}^{x-1} g(k)+N\overline \gl \sum_{k=1}^{x-1} h(k) \right \vert\\
\;\le\; &
\left \vert -N \gl_1 \sum_{k=1}^{x-1} \left[g(k)-h(k) \right] \right \vert
+\left \vert N\left( \overline \gl- \gl_1\right) \sum_{k=1}^{x-1}h(k) \right \vert \\
\;\le\; &
N \gl_1   \sum_{k=1}^{x-1} \left \vert \left[g(k)-h(k) \right] \right \vert  + N  \vert \overline \gl-\gl_1  \vert (x-1)\,.
\end{aligned}
\end{equation}
Then by \eqref{autofuns:approx} and \eqref{gap:asym} with $\gl_1=\gap_N$, we obtain \eqref{der:approx}.

\end{proof}

\subsection{ The shape and derivative in the second segment}

\begin{lemma}\label{lema:gradient}
 If \eqref{LLN} holds, we have
\begin{equation}\label{autofuns:2seg}
\lim_{N \to \infty} \; \max_{x \in \lint \tau_1, \, \tau_2\rint}\;  \left \vert   g(x)- h(x) \right\vert \;=\;0\,,
\end{equation}
and 
\begin{equation}\label{der:2seg}
\lim_{N \to \infty} \; \max_{x \in \lint \tau_1,\, \tau_2\rint}\; N \left \vert  (c\nabla g)(x)-(\nabla h)(x) \right\vert \;=\;0\,.
\end{equation}
\end{lemma}

\begin{rem} The reason why we cannot apply the method of the proof in Lemma \ref{lema:1stapproxfun} is that $B$ is not bounded in the interval $\lint \tau_1, \, \tau_2\rint$.
\end{rem}

\begin{proof}
We first treat \eqref{der:2seg}.
 Recalling $A(x)=1/B(x)$ and $\overline A(x)=1/\overline B(x)$, by \eqref{rec:randeigfun1st} we have
\begin{equation*}
\begin{aligned}
A(x) \left[1-\frac{g(x)}{g(x-1)} \right]&\;=\;r(x-1,x)N^{-1}\,,\\
\overline A(x) \left[1-\frac{h(x)}{h(x-1)} \right]&\;=\;N^{-1}\,.
\end{aligned}
\end{equation*}
Taking the difference of the two equalities above,  we obtain
\begin{equation}\label{Aratiog}
\left \vert  A(x)\frac{g(x)}{g(x-1)}-\overline A(x) \frac{h(x)}{h(x-1)} \right\vert \;\le\; \frac{\vert r(x-1,x)-1 \vert}{N}
+ \vert A(x)-\overline A(x) \vert\,.
\end{equation}
Furthermore, by \eqref{def:blambda}, $B(x)=N b(\gl,x)$ and $\overline B(x)=N \overline b(\gl,x)$ we have
\begin{equation}\label{Angrad}
\begin{aligned}
A(x) \frac{g(x)}{g(x-1)}&\;=\; -\frac{1}{N} \frac{g(x)}{c(x-1,x)\left[ g(x)-g(x-1) \right]}\,,\\
\overline A(x) \frac{h(x)}{h(x-1)}&\;=\;-\frac{1}{N} \frac{h(x)}{\left[ h(x)-h(x-1) \right]}\,.\\
\end{aligned}
\end{equation}
Note that for $x \in \lint \tau_1,  \tau_2 \rint$, by \eqref{fun:h} and the mean value theorem, there exists a constant $M>\pi$ such that
\begin{equation*}
\left\vert \overline A(x) \frac{h(x)}{h(x-1)} \right\vert \;=\;\left\vert -\frac{h(x)}{N \left[ h(x)-h(x-1)\right]} \right\vert \le M-1\,.
\end{equation*}
By \eqref{Aratiog}, \eqref{def:delta1} and Lemma \ref{lema:2seghomA},  for $x \in \lint \tau_1, \tau_2 \rint$,   we also have 
\begin{equation}\label{bound:Ag}
\left\vert  A(x) \frac{g(x)}{g(x-1)} \right\vert \;\le\; M\,.
\end{equation}
 Writing $\Upsilon(x) \,\colonequals\, N (c \nabla g)(x)- N(\nabla h)(x)$, by \eqref{Angrad} and triangle inequality  we have
\begin{equation}
\begin{aligned}\label{diff:2eigfuns}
\vert h(x)-g(x) \vert& \;=\; \left \vert A(x)\frac{g(x)}{g(x-1)} N (c \nabla g)(x)-\overline A(x)\frac{h(x)}{h(x-1) }N (\nabla h)(x) \right \vert\\
&\;\le\;
\left \vert A(x)\frac{g(x)}{g(x-1)} \right\vert \cdot \vert \Upsilon(x) \vert
+ \left\vert A(x)\frac{g(x)}{g(x-1)}- \overline A(x)\frac{h(x)}{h(x-1) } \right \vert  \cdot  \left \vert N (\nabla h)(x) \right\vert\\
&\;\le\; M  \vert \Upsilon (x) \vert+ M \tilde \delta_N
\end{aligned}
\end{equation}
where we have used \eqref{bound:Ag}, $\vert N (\nabla h)(x) \vert \le M$, \eqref{Aratiog} and
\begin{equation}\label{def:tilde-delta}
\tilde \delta_N \;\colonequals\; \delta_N^{(0)} + \max_{x \in \lint \tau_1, \, \tau_2 \rint} \vert A(x)-\overline A(x) \vert\,.
\end{equation}
Note that
$\tilde \delta_N$ tends to zero as $N \to \infty$ by \eqref{LLN} and Lemma \ref{lema:2seghomA}. In order to show that the r.h.s.\ of \eqref{diff:2eigfuns} tends to zero, we are left to  show that $\Upsilon(x)$ tends to zero. In view of \eqref{grad:g}, we define 
\begin{equation}\label{def:sum-difference}
s(x) \;\colonequals\; \sum_{k= 1}^x \left\vert g(k)- h(k) \right\vert\,, \quad \forall \, x \in \lint 1, N+1\rint\,.
\end{equation} 
By \eqref{diff:2eigfuns}, \eqref{diff:derivatives} and \eqref{grad:g}, for $x \in \lint \tau_1\,,\tau_2\rint$  we have
\begin{equation*}
\begin{aligned}
s(x) &\;=\; s(x-1)+ \vert h(x)-g(x) \vert\\
&\;\le\; s(x-1)+ M \vert  \Upsilon(x)\vert+ M \tilde \delta_N \\
&\;\le\; s(x-1)+ M \left[ N \gl_1 s(x-1)+N \vert \overline \gl-\gl_1 \vert (x-1) \right]+M \tilde \delta_N\\
&\;=\;\left(1+M N \gl_1 \right) s(x-1)+M N \vert \overline \gl-\gl_1\vert (x-1)+ M \tilde \delta_N\,.
\end{aligned}
\end{equation*}
Iterating the inequality above, for $x \in \lint \tau_1,\, \tau_2 \rint$ we obtain
\begin{equation}
\begin{aligned}\label{rec:sx}
s(x) & \;\le\; \left( 1+M N \gl_1\right)^N s(\tau_1)+ \left( 1+M N \gl_1\right)^N \sum_{k=\tau_1}^{x-1} \left[ M N \vert \overline \gl-\gl_1 \vert k+ M\tilde \delta_N \right]\\
&
\;\le\; C(M, \ga) \left[ s(\tau_1)+  MN^3 \vert \overline \gl -\gl_1\vert +  M N \tilde \delta_N \right] \\
&
\;\le\;
C(M, \ga) \left[ N \max_{x \in \lint 1,\, \tau_1 \rint} \vert g(x)-h(x) \vert+  MN^3 \vert \overline \gl -\gl_1\vert +  M N \tilde \delta_N\right]\,,
\end{aligned}
\end{equation}
where we have used $(1+MN \gl_1)^N \le C(M, \ga)$ and the last term above does not depend on $x$. 
Plugging \eqref{rec:sx} into \eqref{diff:derivatives} which holds for all $x \in \lint 1,\, N+1\rint$,  we obtain
\begin{equation*}
\begin{aligned}
\max_{x \in \lint \tau_1,\, \tau_2 \rint} \; \vert N (c \nabla g)(x)-N(&\nabla h)(x) \vert
\;\le\:  N^2 \vert \overline \gl-\gl_1 \vert \\
&+
N \gl_1  C(M, \ga)\left[ N \max_{x \in \lint 1,\, \tau_1 \rint} \vert g(x)-h(x) \vert+  MN^3 \vert \overline \gl -\gl_1\vert + M N \tilde \delta_N \right]
\end{aligned}
\end{equation*}
which tends to zero as $N \to \infty$ by  Lemma \ref{lema:1stapproxfun}, \eqref{gap:asym} and \eqref{def:tilde-delta}. Therefore we prove \eqref{der:2seg}.
Then by \eqref{diff:2eigfuns} and \eqref{der:2seg}, we conclude the proof of \eqref{autofuns:2seg}.

\end{proof}

\subsection{The shape and derivative  in the last segment}  As $B$ is bounded in the segment $\lint \tau_2, N+1 \rint$, we just need to apply \eqref{iterate:ux} to repeat exactly the same proof for \eqref{autofuns:approx}. Concerning the approximation of the two derivatives, by  \eqref{diff:derivatives} and \eqref{eigfun:shape} we conclude the proof.

\begin{proof}[Proof of \eqref{eigfun:shape} and \eqref{approx:der2eigfuns} in Theorem \ref{th:gapshapeder}]
We combine Lemma \ref{lema:1stapproxfun}, Lemma \ref{lema:gradient} and the comments above about the last segment $\lint \tau_2,\, N \rint$ to conclude the proof of \eqref{eigfun:shape} and \eqref{approx:der2eigfuns}.

\end{proof}

\bibliographystyle{alpha}
\bibliography{library.bib}

\end{document}